\documentclass{elsart}
\usepackage{amssymb}
\journal{Journal of Algebra}
\begin{document}
\begin{frontmatter}
\title{On open normal subgroups of parahorics}
\author{Fran\c cois Court\`es} 
\address{Universit\'e de Poitiers\\
D\'epartement de Math\'ematiques\\
UMR 6086 du CNRS\\
T\'el\'eport 2\\
Boulevard Marie et Pierre Curie\\
86962 Futuroscope Chasseneuil Cedex}
%\ead{courtes@math.univ-poitiers.fr}
\abstract Let $F$ be a local complete field with discrete valuation, and let $G$ be a quasi-split group over $F$ which splits over some unramified extension of $F$. Let $P$ be a parahoric subgroup of the group $G(F)$ of $F$-points of $G$; the open normal pro-nilpotent subgroups of $P$ can be classified using the standard normal filtration subgroups of Prasad and Raghanathan. More precisely, we show that if $G$ is quasi-simple and satisfies some additional conditions, $H$ is, modulo a subgroup of some maximal torus of $G$, either one of these filtration subgroups or the product of one of them by a standard normal filtration subgroup of $P\cap M$, where $M$ is a proper Levi subgroup of $G$.
\endabstract
\keyword Reductive groups over local fields \sep parahoric subgroups \sep normal subgroups
\MSC 20G25 \sep 20E07
\endkeyword
\end{frontmatter}
\def \mth{\mathbb}
\newtheorem{theo}{Theorem}[section]
\newtheorem{propo}[theo]{Proposition}
\newtheorem{lemme}[theo]{Lemma}
\newtheorem{coro}[theo]{Corollary}

\section{Introduction}

Let $F$ be a local complete field with discrete valuation; let $\mathcal{O}$ be its ring of integers, ${\mathfrak{p}}$ the maximal ideal of $\mathcal{O}$, $K=\mathcal{O}/{\mathfrak{p}}$ its residual field. Let $\varpi$ be an uniformizer of $F$, and $v_F$ be the normalized valuation on $F$.

Let $G$ be a connected reductive algebraic group defined over $F$; we'll assume $G$ is quasisplit, and splits over some unramified extension of $F$; let $G(F)$ be the set of $F$-points of $G$. Let $\Phi$ be the root system of $G(F)$ relatively to some maximal split torus $S$ of $G$; we'll assume that $\Phi$ satisfies the following condition:

{\em (G1): if the characteristic $p$ of $K$ is $2$, $\Phi$ is simply-laced; if $p=3$, $\Phi$ has no connected component of type $G_2$.}

Let $G^0$ be the parahoric component of $G(F)$, that is the subgroup of $G(F)$ generated by all parahoric subgroups of $G(F)$. $G^0$ is an open normal subgroup of $G(F)$, and the quotient $G(F)/G^0$ is abelian; moreover, $G(F)/G^0$ is finite if and only if $G$ is semisimple, and $G^0=G(F)$ if and only if $G$ is simply-connected.

Let $\mathcal{B}$ be the Bruhat-Tits building of $G(F)$, and let $\mathcal{A}$ be the apartment of $\mathcal{B}$ associated to $S$; $\mathcal{A}$ is isomorphic as a ${\mth{R}}$-affine space to $(X_*(S)\otimes{\mth{R}})/(X_*(Z)\otimes{\mth{R}})$, where $Z$ is the split component of the center of $G$ and $X_*(S),X_*(Z)$ are the groups of cocharacters of $S,Z$. The subgroup $Q$ of the group $X^*(S)$ of characters of $S$ generated by $\Phi$ can then be viewed as a group of affine functions on $\mathcal{A}$, by the standard duality product.

For every $x\in{\mth{B}}$ and every $r\in{\mth{R}}^+$, let $B(x,r)$ be the subset of elements $y$ of $\mathcal{B}$ satisfying the following property: for every apartment $\mathcal{A}'$ containing both $x$ and $y$ and every root $\alpha$ of $G$ relatively to $S'$, $S'$ being the maximal split torus associated to $\mathcal{A}'$, we have $\alpha(y)-\alpha(x)\leq r$. Since $B(x,r)$ is nonempty, the subgroup $G_{x,r}$ of elements of $G^0$ fixing $B(x,r)$ pointwise is an open bounded subgroup of $G(F)$. Moreover, if $G_x$ is the parahoric subgroup of $G(F)$ fixing $x$, $G_{x,r}$ is normal in $G_x$, and $G_{x,s}\subset G_{x,r}$ for $s>r$. We'll also write:
\[G_{x,r}^+=\bigcup_{s>r}G_{x,s};\]
the group $G_{x,r}^+$ is also normal in $G_x$. Note that since the valuation is discrete, we have $G_{x,r}^+=G_{x,s}$ for any $s>r$ sufficiently close to $r$.

These groups were first introduced in \cite{pr} and have been used in \cite{mp1} and \cite{mp2} to classify unramified types, in the context of $p$-adic fields, that is fields with finite residual fields. They constitute a standard filtration of parahoric subgroups of $G(F)$ by open normal subgroups, in the following sense: let $A$ be a facet of $\mathcal{B}$, and $P_A$ be the parahoric subgroup of $G$ fixing $A$ pointwise. Then for every $x\in A$, $G_{x,0}=G_x=P_A$, and the $G_{x,r}$, with $x\in A$ and $r\geq 0$, are a basis of neighborhoods of unity. Moreover, for every $x\in A$, $G_{x,0}/G_{x,0}^+$ is the group of $K$-points of a reductive $K$-group (which depends only on $A$), and for every $r>0$, $G_{x,r}/G_{x,r}^+$ is abelian.

We can consider the concave function $f$ on $\Phi$ defined by $f(\alpha)=\alpha(x)+r$ for every $\alpha\in\Phi$. Let $U_{x,r}=U_{G,x,r}$ be the subgroup $U_f$ of $G(F)$ attached to $f$ as in \cite[I.6.4]{bt}; $U_{x,r}$ is a normal subgroup of $G_x$, and there exists a subgroup $T'$ of $T(F)$, where $T$ is the centralizer of $S$ in $G$, such that $G_{x,r}=T'U_{x,r}$. (For given $x,r$, the subgroups $T''U_{x,r}$, with $T''\subset T(F)$, which are normal in $G_x$ may be classified with the help of \cite[I.6.4.19]{bt}.)

Ona can naturally ask if these subgroups, for $r>0$, exhaust the open normal pro-nilpotent subgroups of $G_x$, at least when $G$ is quasi-simple. The answer is in general no, but they can still be used to classify them.

First we'll show that to every open bounded subgroup $H$ of $G$ normalized by the unique parahoric subgroup $P_T$ of $T(F)$, we can attach a concave function $f_H$ on $\Phi$. More precisely, we have the following result:

\begin{propo}\label{nsgf}
Assume $\Phi$ satisfies (G1), and:
\begin{itemize}
\item either $p\neq 2$ or $F$ is absolutely unramified;
\item $K$ has at least $4$ elements, and if $K$ has exactly $4$ (resp. $5$) elements, $\Phi$ has no component of type $A_2$ (resp. $A_1$, $C_n$ or $BC_n$).
\end{itemize}
Let $H$ be an open bounded pro-nilpotent subgroup of $G(F)$ normalized by $P_T$; there exists a bounded subgroup $T'$ of $T(F)$ and a concave function $f_H$ on $\Phi$ such that $H=T'U_{f_H}$.
\end{propo}

Since every parahoric subgroup of $G(F)$ contains the parahoric subgroup of some maximal torus of $G(F)$, this will in particular be true for open normal subgroups of $G_x$.

We then consider more precisely the concave functions attached to the open normal pro-nilpotent subgroups of parahorics; and more particularly of some given Iwahori subgroup $I$ of $G(F)$; we will in fact study the function $f_\varepsilon=f-f_I$. By establishing some properties of the elements of a given class of subsets of $\Phi$, which will be called $\Delta'$-complete subsets, where $\Delta'$ is the extended basis of $\Phi$ associated to $I$, we will then determine $r\in{\mth{R}}_+^*$ and $x\in\mathcal{A}$ which will satisfy the main result of the paper:

\begin{theo}\label{cln}
Assume $G$ is quasi-simple and the conditions of the previous proposition are satisfied. Let $A$ be a facet of $\mathcal{B}$, $\overline{A}$ be the closure of $A$ and $G_A$ be the parahoric subgroup of $G$ fixing $A$; let $H$ be any open normal pro-nilpotent subgroup of $G_A$. There exists $x\in\overline{A}$ and $r\in{\mth{R}}^*_+$ such that $H$ satisfies one of the following conditions:
\begin{itemize}
\item $U_{x,r}\subset H\subset G_{x,r}$;
\item there exists $r'>r$ and a proper Levi subgroup $M$ of $G$ containing $T(F)$ and such that $U_{x,r'}U_{M,x,r}\subset H\subset G_{x,r'}M_{x,r}$.
\end{itemize}
\end{theo}

We conclude with an easy generalization of our main result to non-quasi-simple groups.

\section{Generalities}

\subsection{A few miscellaneous notations}

Let $F'$ be a valued field with discrete valuation; we denote by $v_{F'}$ the normalized valuation on $F'$, that is the unique valuation whose image is precisely ${\mth{Z}}\cup\{+\infty\}$. If $E$ is any extension of $F'$, we'll denote again by $v_{F'}$ the valuation on $E$ whose restriction to $F'$ is $v_{F'}$.

Let $H$ be a group and $H'$ be a subgroup of $H$; we'll write $N_H(H')$ (resp. $Z_H(H'))$ for the normalizer (resp. the centralizer) of $H'$ in $H$.

For every $r\in{\mth{R}}$, we'll write $\operatorname{floor}(r)$ (resp. $\operatorname{ceil}(r)$) for the largest integer $\leq r$ (resp. the smallest integer $\geq r$).

\subsection{More about roots and affine roots}

This section is devoted to some results about root systems which will be useful later. (For the basic results about root systems, see \cite{bou}).

Let $V$ be a finite-dimensional ${\mth{R}}$-vector space, and let $\Phi$ be a root system contained in $V$. We'll assume in this subsection that $\Phi$ is irreducible.

We'll denote by $W$ the Weyl group of $\Phi$ . Let $(.,.)$ be a $W$-invariant scalar product on $V$; we'll assume $(.,.)$ is chosen such that:
\begin{itemize}
\item if $\Phi$ is simply-laced, $(\alpha,\alpha)=2$ for every $\alpha\in\Phi$;
\item if $\Phi$ is not simply-laced, $(\alpha,\alpha)=2$ if $\alpha$ is a long root in $\Phi$.
\end{itemize}

If $\Phi$ is reduced, we'll write $h(\Phi)$ for the Coxeter number of $\Phi$. If $\Phi$ is of type $BC_n$, we'll set $h(\Phi)=2n+1$, that is the Coxeter number of $\Phi$ plus one.

Let $\Phi_{aff}$ be the affine root system associated to $\Phi$, which will be identified with the subset $\Phi\times{\mth{Z}}$ of $X^*(S)\times{\mth{Z}}$. We'll denote by $W_{aff}$ the affine Weyl group of $\Phi_{aff}$.

Let $\Delta$ be a basis of $\Phi$, and let $\Delta'$ be the extended basis $\Delta\cup\{-\alpha_M\}$ of $\Phi$, where $\alpha_M$ is the largest root of $\Phi$ w.r.t $\Delta$. Set:
\[\Delta_{aff}=\{(\alpha,0)|\alpha\in\Delta\}\cup\{(-\alpha_M,1)\};\]
$\Delta_{aff}$ is a basis of $\Phi_{aff}$.

Let $\Phi^+$ (resp. $\Phi^+_{aff}$) be the set of positive elements of $\Phi$ (resp. $\Phi_{aff}$) w.r.t $\Delta$ (resp. $\Delta_{aff}$); $\Phi^+_{aff}$ is the set of elements $(\alpha,v)$ in $\Phi_{aff}$ such that $v\geq 0$ if $\alpha\in\Phi^+$ and $v\geq 1$ else. It is well-known (see for example \cite[2]{pr}) that every element of $\Phi^+_{aff}$ can be uniquely written as a linear combination with nonnegative integer coefficients of elements of $\Delta_{aff}$; we deduce from this the following result:

\begin{propo}
Every element $\alpha$ of $\Phi$ can be written as a linear combination with nonnegative integer coefficients of elements of $\Delta'$; moreover, there is an unique such combination satisfying the following condition: the sum of the coefficients is strictly smaller than $h(\Phi)$.
\end{propo}

Consider the projection on $\Phi$ of the decomposition of $(\alpha,v)\in\Phi^+_{aff}$ for any suitable $v$; the unique combination satisfying the last condition is the one obtained with $v=0$ (resp. $v=1$) if $\alpha>0$ (resp. $\alpha<0$). $\Box$

We'll call the height of $(\alpha,v)$ (resp. $\alpha$) relatively to $\Delta_{aff}$ (resp. $\Delta'$) and we'll write $h(\alpha,v)$ (resp. $h(\alpha)$) for the sum of the coefficients of $(\alpha,v)$ (resp. $\alpha$) defined as above.

For every $\alpha\in\Phi$, let $\varepsilon_\alpha$ be the smallest integer such that $(\alpha,\varepsilon_\alpha)\in\Phi^+_{aff}$, that is $\varepsilon_\alpha=0$ if $\alpha>0$, $\varepsilon_\alpha=1$ if $\alpha<0$. The following results are immediate:

\begin{itemize}
\item if $(\alpha,u),(\beta,v)$ are elements of $\Phi^+_{aff}$ such that $(\alpha+\beta,u+v)\in\Phi_{aff}$, then $h(\alpha+\beta,u+v)=h(\alpha,u)+h(\beta,v)$;
\item for every $\alpha\in\Phi$, $h(-\alpha,\varepsilon_{-\alpha})=h(\Phi)-h(\alpha,\varepsilon_\alpha)$.
\end{itemize}

We obtain the following corresponding assertions about elements of $\Phi$:

\begin{itemize}
\item if $\alpha,\beta$ are elements of $\Phi$ such that $\alpha+\beta\in\Phi$, then $h(\alpha+\beta)\equiv h(\alpha)+h(\beta)$ modulo $h(\Phi)$;
\item for every $\alpha\in\Phi$, $h(-\alpha)=h(\Phi)-h(\alpha)$.
\end{itemize}

In particular, for every $\alpha\in\Phi$, we have $h(\alpha)\in\{1,\dots,h(\Phi)-1\}$.

We can define a partial order on $\Phi^+_{aff}$ relatively to $\Delta_{aff}$ as follows: for every $(\alpha,u),(\beta,v)\in\Phi$, we have $\alpha\leq\beta$ if all the coefficients of elements of $\Delta_{aff}$ in the decomposition of $(\alpha,u)$ are smaller than or equal to the corresponding ones in the decomposition of $(\beta,v)$. In particular, if $u<v$, we have $(\alpha,u+\varepsilon_\alpha)<(\beta,v+\varepsilon_\beta)$ for every $\alpha,\beta\in\Phi$.

We deduce from this order a partial order on $\Phi$ as follows: $\alpha\leq\beta$ if and only if $(\alpha,v+\varepsilon_\alpha)\leq(\beta,v+\varepsilon_\beta)$ for every $v$. This partial order is different from the usual one relative to $\Delta$, but it is easy to see that the restrictions of both orders to the subset of positive (resp. negative) elements are the same; the minimal elements (resp. the maximal elements) of our order are the elements of $\Delta'$ (resp. $-\Delta'$), and we easily see that $\alpha\leq\beta$ implies $-\beta\leq-\alpha$. In the sequel, unless another order on $\Phi$ is explicitly mentioned, we will always refer to this one.

We have the following results:

\begin{propo}
Let $\alpha,\beta$ be two elements of $\Phi$ such that $\alpha+\beta\in\Phi$; if $h(\alpha+\beta)=h(\alpha)+h(\beta)$ (resp. if $h(\alpha+\beta)=h(\alpha)+h(\beta)-h(\Phi)$), then $\alpha+\beta$ is greater (resp. lesser) than both $\alpha$ and $\beta$.
\end{propo}

The first assertion is obvious. The second one is obtained by simply remarking that if $h(\alpha+\beta)=h(\alpha)+h(\beta)-h(\Phi)$, then we have:
\[h(-\alpha-\beta)=h(\Phi)-h(\alpha+\beta)\]
\[=2h(\Phi)-h(\alpha)-h(\beta)=h(-\alpha)+h(-\beta).\]
$\Box$

\begin{propo}\label{decprg}
Let $\alpha,\beta$ be elements of $\Phi$; we have $\alpha\leq\beta$ if and only if there exist elements $\gamma_1,\dots,\gamma_t$ of $\Delta'$ such that:
\begin{itemize}
\item for every $i\in\{1,\dots,t\}$, $\alpha_i=\alpha+\sum_{j=1}^i\gamma_i$ is an element of $\Phi$;
\item $\alpha_t=\beta$.
\end{itemize}
\end{propo}

Assume $\alpha\leq\beta$; we'll show the existence of the $\gamma_i$. Write $\beta=\alpha+\sum_{i=1}^t\delta_i$, with $\delta_1,\dots,\delta_t\in\Delta'$. If $t=0$, there is nothing to prove; suppose $t>0$. Since $(\beta,\beta)>0$, we have either $(\beta,\delta_i)>0$ for some $i$ or $(\beta,\alpha)>0$. In the first case, set $\gamma_t=\delta_i$; the result follows from the induction hypothesis applied to $\alpha$ and $\beta-\delta_i$. In the second case, $\beta-\alpha$ is a root. If $(\alpha,\delta_i)<0$ for some $i$, we can set $\gamma_1=\delta_i$ and conclude by the induction hypothesis applied to $\alpha+\delta_i$ and $\beta$; if there is no such $i$, we are in one of the following two cases:
\begin{itemize}
\item $\Phi$ is of type $B_n$, $C_n$ or $F_4$, and $\alpha$ and $\beta-\alpha$ are orthogonal but not strongly orthogonal;
\item $\Phi$ is of type $G_2$, and $\alpha$ and $\beta-\alpha$ are short roots such that $(\alpha,\beta)=1$.
\end{itemize}
In both cases, $\beta$ is a long root; the result is then true for $-\beta$ and $-\alpha$, which obviously implies the result for $\alpha$ and $\beta$.

Conversely, assume the $\gamma_i$ do exist. By an easy induction we may assume $t=1$; we then have either $h(\beta)=h(\alpha)+1$ or $h(\beta)=h(\alpha)-h(\Phi)+1$; since the length of any element of $\Phi$ is contained in $\{1,\dots,h(\Phi)-1\}$, the second case is impossible and we obtain $\alpha\leq\beta$, as required. $\Box$

Let $\overline{\Phi}$ be the subset $\Phi\cup\{0\}$ of $X^*(T)$. Let $\alpha,\beta$ be two elements of $\overline{\Phi}$; it is easy, with the help of \cite[1.3, cor. to th. 1]{bou}, to check that we always have:
\begin{itemize}
\item if $(\alpha,\beta)<0$, then $\alpha+\beta\in\overline{\Phi}$;
\item if $(\alpha,\beta)>0$, then $\alpha-\beta\in\overline{\Phi}$.
\end{itemize}

Set also $\overline{\Phi_{aff}}=\overline{\Phi}\times{\mth{Z}}$. For every $v>0$, $(0,v)\in\overline{\Phi_{aff}}$ can also be uniquely written as a linear combination with nonnegative integer coefficients of elements of $\Delta_{aff}$; moreover, we have $h(0,v)=h(\Phi)v$, and the partial order on $\Phi^+_{aff}$ extends canonically to $\overline{\Phi^+_{aff}}=\Phi^+_{aff}\cup(\{0\}\times{\mth{N}}^*)$ by setting $(\alpha,v-1+\varepsilon_\alpha)\leq(0,v)\leq(\beta,v+\varepsilon_\beta)$ for every $\alpha,\beta$.

\begin{propo}\label{infsum}
Let $(\alpha_1,c_1+\varepsilon_{\alpha_1}),\dots,(\alpha_t,c_t+\varepsilon_{\alpha_t})$ be elements of $\Phi^+_{aff}$ whose sum is also an element $(\alpha,c+\varepsilon_\alpha)$ of $\Phi^+_{aff}$. Let $\beta\in\overline{\Phi}$ and $c'\in\{c,c+1\}$ be such that $(\alpha,c+\varepsilon_\alpha)\leq(\beta,c'+\varepsilon_\beta)\leq(0,c+1)$ (resp. $(0,c)\leq(\beta,c'+\varepsilon_\beta)\leq(\alpha,c+\varepsilon_\alpha)$); there exist then elements $(\beta_1,c'_1+\varepsilon_{\beta_1})\dots,(\beta_t,c'_t+\varepsilon_{\beta_t})$ of $\overline{\Phi^+_{aff}}$ whose sum is $(\beta,c+\varepsilon_\beta)$ and such that $(\alpha_i,c_i+\varepsilon_{\alpha_i})\leq(\beta_i,c'_i+\varepsilon_{\beta_i})\leq(0,c_i+1)$ (resp. $(0,c_i)\leq(\beta_i,c'_i+\varepsilon_{\beta_i})\leq(\alpha_i,c_i+\varepsilon_{\alpha_i})$) for every $i$.
\end{propo}

We'll show the proposition with $(\alpha,c+\varepsilon_\alpha)\leq(\beta,c'+\varepsilon_\beta)\leq(0,c+1)$, the proof for $(0,c)\leq(\beta,c'+\varepsilon_\beta)\leq(\alpha,c+\varepsilon_\alpha)$ being similar. With the previous proposition and an easy induction, we can assume there exists $\gamma\in\Delta'$ such that $\beta=\alpha+\gamma$, hence $(\beta,c'+\varepsilon_\beta)=(\alpha,c+\varepsilon_\alpha)+(\gamma,\varepsilon_\gamma)$. If $(\alpha_i,\gamma)<0$ for some $i$, then $\alpha_i+\gamma\in\Phi$, and since $(\beta_i,c'_i+\varepsilon_{\beta_i})=(\alpha_i,c_i+\varepsilon_{\alpha_i})+(\gamma,\varepsilon_\gamma)>(\alpha_i,c_i+\varepsilon_{\alpha_i})$, setting $(\beta_j,c'_j+\varepsilon_{\beta_j})=(\alpha_j,c_j+\varepsilon_{\alpha_j})$ for every $j\neq i$ proves the proposition. Suppose now $(\alpha_i,\gamma)\geq 0$ for every $i$. Since $(\alpha,\alpha)>0$, there exists an $i$ such that $(\alpha,\alpha_i)>0$, hence $\alpha-\alpha_i\in\overline{\Phi}$. If $\alpha=\alpha_i$, then $\alpha_i+\gamma\in\Phi$ and we conclude as before; we may then suppose $\alpha-\alpha_i\neq 0$.

Since $(\alpha,\alpha_i)>0$ and $(\gamma,\alpha_i)\geq 0$, we have $(\beta,\alpha_i)=(\alpha+\gamma,\alpha_i)>0$, hence $\beta-\alpha_i\in\overline{\Phi}$. Moreover, setting $(\alpha',d+\varepsilon_{\alpha'})=(\alpha,c+\varepsilon_\alpha)-(\alpha_i,c_i+\varepsilon_{\alpha_i})$ and $(\beta',d'+\varepsilon_{\beta'})=(\beta,c'+\varepsilon_\beta)-(\alpha_i,c_i+\varepsilon_{\alpha_i})$, we must have either $d=d'$ or $\beta'=0$, since $d<d'$ and $\beta'\neq 0$ would imply $h(\beta',d'+\varepsilon_{\beta'})\geq h(\alpha',d+\varepsilon_{\alpha'})+2$, hence $h(\beta,c+\varepsilon_\beta)\geq h(\alpha,c+\varepsilon_\alpha)+2$, which contradicts our assumptions; we then have $(\alpha',d+\varepsilon_{\alpha'})\leq(\beta',d'+\varepsilon_{\beta'})\leq(0,d+1)$. This way, we can show the result by induction on $t$; since the case $t=1$ is trivial, the proposition is proved. $\Box$

\subsection{Valuations}

From now on, we'll assume $\Phi$ is the root system of $G$ relatively to some maximal split torus $S$ of $G$. We will denote by $T$ the centralizer of $S$ in $G$; since $G$ is quasisplit, $T$ is a maximal torus of $G$, which splits over the same unramified extension of $F$ as $G$.

For every $\alpha\in\Phi$, let $U_\alpha$ be the root subgroup of $G$ associated to $\alpha$. If $2\alpha$ is not a root, there exists a finite unramified extension $F_\alpha$ of $F$ such that the group $U_\alpha(F)$ is isomorphic to $F_\alpha$; for every $x\in F_\alpha$, we'll write $u_\alpha(x)$ for the image of $x$ in $U_\alpha(F)$ by sone given isomorphism. If $2\alpha\in\Phi$, there exists a quadratic unramified extension $F_\alpha$ of $F_{2\alpha}$ such that $U_\alpha(F)$ is isomorphic to the semidirect product $H_\alpha$ of $F_\alpha$ by $F_{2\alpha}$, with the addition law defined as follows:
\[(x,y)+(x',y')=(x+x',y+y'+x\sigma(x)),\]
when $\sigma$ is the nontrivial element of $\operatorname{Gal}(F_\alpha/F_{2\alpha})$; moreover, the image of $\{(0,y)|y\in F_{2\alpha}\}$ is $U_{2\alpha}(F)$. Let's choose such an isomorphism; for every $x\in F_\alpha$, we'll set $u_\alpha(x)$ to be the image of $(x,0)$ by that isomorphism.

For every $\alpha$, we'll write $\mathcal{O}_\alpha$ for the ring of integers of $F_\alpha$, ${\mathfrak{p}}_\alpha$ for the maximal ideal of $\mathcal{O}_\alpha$, $K_\alpha$ for the residual field of $F_\alpha$.

Let $(v_\alpha)_{\alpha\in\Phi}$ be a valuation on $(T,(U_\alpha)_{\alpha\in\Phi})$; for every $\alpha$, $v_\alpha$ is an application from $U_{\alpha}(F)$ to ${\mth{R}}\cup\{+\infty\}$ satisfying the following conditions:

\begin{itemize}
\item for every $\alpha\in\Phi$, there exist constants $b_\alpha\in{\mth{R}}^*_+$, $c_\alpha\in{\mth{R}}$ such that for every $x\in F_\alpha$, $v_\alpha(u_\alpha(x))=b_\alpha v_F(x)+c_\alpha$;
\item if $2\alpha\in\Phi$, for every $(x,y)\in H_\alpha$, if $u$ is the image of $(x,y)$ in $U_\alpha(F)$, $v_\alpha(u)=\operatorname{Inf}(v_\alpha(x,0),\frac 12v_{2\alpha}(0,y))$;
\item the commutator relations: for every $\alpha\in\Phi$ and every $r\in{\mth{R}}$, let $U_{\alpha,r}$ be the subgroup of elements $u$ of $U_\alpha(F)$ such that $v_\alpha(u)\geq r$. If $\alpha,\beta$ are elements of $\Phi$ such that $\alpha+\beta\in\Phi$, and if $r,r'$ are elements of ${\mth{R}}$, the commutator subgroup $[U_{\alpha,r},U_{\beta,r'}]$ is contained in $U_{\alpha+\beta,r+r'}$; (and by \cite{chev} and \cite[Appendix A]{bt}, if (G1) is satisfied, this inclusion is an equality when $r\in \operatorname{Im}(v_\alpha)$ and $r'\in \operatorname{Im}(V_\beta)$);
\item for every $\alpha\in\Phi$, if $u\in U_\alpha(F)$ and $u'\in U_{-\alpha}(F)$ generate a subgroup of $G$ which is bounded but not pro-solvable, then $v_\alpha(u)+v_{-\alpha}(u')=0$.
\end{itemize}

This is a rewriting of \cite[I, definition 6.2.1]{bt} in our context (with a discrete valuation and a group which splits over an unramified extension).

We will also define a valuation on $T(F)$ the following way: for each $r\in{\mth{R}}$, we define the group $T_r$ as the subgroup of elements $t\in T(F)$ such that for every $\alpha\in\Phi$ and every $r'\in{\mth{R}}$, we have:
\[[t,U_{\alpha,r'}]\subset U_{\alpha,r+r'}.\]
For every $t\in T(F)$, we'll set $v_0(t)=r$ if $t$ belongs to $T_r$ but not to any $T_{r'}$, $r'>r$. It is easy to see that adding $v_0$ to the family $(v_\alpha)_{\alpha\in\Phi}$ leads to an extension of the valuation in the sense of \cite[I.6.4.38]{bt}.

We can even define valuations on bounded subgroups of $G(F)$. For every $\alpha\in\Phi$, set $U'_{\alpha}=U_\alpha(F)$ if $2\alpha\not\in\Phi$, and if $2\alpha\in\Phi$, let $U'_\alpha$ be the image of the application $u_\alpha$ (which is not necessarily a subgroup of $U_\alpha(F)$); set $U'_{\alpha,r}=U_{\alpha,r}\cap U'_\alpha$ for every $r\in{\mth{R}}^+$. Let $G(F)_r$ be, for a given $r$, the subgroup of $G(F)$ generated by $T_r$ and the $U'_{\alpha,r}$, $\alpha\in\Phi$; we'll also write $G(F)_r^+$ for the union of all the $G(F)_{r'}$, $r'>r$. The group $G(F)_0$ is a parahoric subgroup of $G(F)$, and $G(F)_0^+$ is its pro-nilpotent radical; moreover, the groups $G(F)_r$ and $G(F)_r^+$ are normal subgroups of $G$.

Of course, these groups depend on the chosen valuation. In fact, we can even show that we can attach to any valuation a point $x$ of the apartment $\mathcal{A}$ of $\mathcal{B}$ associated to $T$, and conversely, in such a way that $G(F)_r=G_{x,r}$ for every $r$, but we will not use this fact; we will only use the following property (see \cite[I.7.2]{bt}): for every parahoric subgroup $H$ of $G(F)$ containing the parahoric subgroup $P_T$ of $T(F)$, there exists a valuation such that $H=G(F)_0$, which implies by \cite[I.3.3.1]{bt} that every bounded subgroup of $G(F)$ is contained in the group $G(F)_0$ relative to some valuation.

For every $g\in G(F)_0$, we'll set $v_G(g)=r$ if $g$ belongs to $G(F)_r$ but not to $G(F)_r^+$.

For convenience, we'll write $U_0=T$, and $U_{0,r}=T_r$ for every $r$. The following result follows immediately from the definitions:

\begin{propo}
Let $\alpha$ be any element of $\overline{\Phi}$, and let $u$ be any element of $U'_\alpha\cap G(F)_0$. Then $v_G(u)=v_\alpha(u)$.
\end{propo}

\subsection{Concave functions}

From now on and until the end of the paper, we'll assume $\Phi$ satisfies (G1). Moreover, in this subsection, we'll also assume $\Phi$ is connected.

Let $f$ be a map from $\Phi$ to ${\mth{R}}$; $f$ is said to be concave if it satisfies the following conditions:

\begin{itemize}
\item for each $\alpha,\beta\in\Phi$ such that $\alpha+\beta\in\Phi$, $f(\alpha+\beta)\leq f(\alpha)+f(\beta)$;
\item for every $\alpha\in\Phi$, $f(\alpha)+f(-\alpha)\geq 0$.
\end{itemize}

Let $U_f$ be the subgroup of $G$ generated by the $U_{\alpha,f(\alpha)}$, $\alpha\in\Phi$. The following properties are well-known (see \cite[I.6.4]{bt}):

\begin{itemize}
\item if $f$ is concave, then for every $\alpha\in\Phi$, $U_f\cap U_\alpha(F)$ is equal to $U_{\alpha,f(\alpha)}$ if $2\alpha\not\in\Phi$, to $U_{\alpha,f(\alpha)}U_{2\alpha,f(2\alpha)}$ if $2\alpha\in\Phi$;
\item since $\Phi$ satisfies (G1), the converse is true if $f$ is optimal (i.e. if for every $\alpha\in\Phi$, $f(\alpha)\in \operatorname{Im}(v_\alpha)$).
\end{itemize}

Moreover, if $f$ is concave and for every $\alpha\in\Phi$, $f(\alpha)+f(-\alpha)>0$, we deduce immediately from \cite[I.6.4.10]{bt} that $U_f$ is pro-solvable and that we have:
\[U_f=(T\cap U_f)\prod_{\alpha\in\Phi}(U_\alpha(F)\cap U_f)\]
the product being taken in any order. Conversely, for any pro-solvable subgroup $H$ of $G$ satisfying the above condition and such that for every $\alpha\in\Phi$, $U_\alpha(F)\cap H=U_{\alpha,f_\alpha}$ for some $f_\alpha$, if moreover $f_\alpha$ is maximal among the elements of ${\mth{R}}$ satisfying that property, then the map $f_H:\alpha\mapsto f_\alpha$ is concave and optimal, we have $f_H(\alpha)+f_H(-\alpha)>0$ for every $\alpha\in\Phi$, and there exists a subgroup $T'$ of $T(F)$ such that $H=T'U_{f_H}$.

Assume now the valuation has been chosen in such a way that for every $\alpha$, $v_\alpha(U'_\alpha)=c_\alpha+{\mth{Z}}\cup\{+\infty\}$ for some constant $c_\alpha$; since $G$ splits over an unramified extension of $F$, this is always possible.

\begin{propo}\label{ctpar}
Let $f$ be a concave function; the group $P_f=P_TU_f$ is contained in a maximal parahoric subgroup. Moreover, if $f$ is optimal and for every $\alpha\in\Phi$, we have $f(\alpha)+f(-\alpha)\leq 1$ (resp. $f(\alpha)+f(-\alpha)=1$), $P_f$ is a parahoric subgroup (resp. an Iwahori subgroup) of $G(F)$.
\end{propo}

According to \cite[I.6.4.9]{bt}, $P_f$ is the finite union of the cosets:
\[P_T(\prod_{\alpha\in\Phi}U_{\alpha,f(\alpha)})n,\]
where $n$ belongs to a system of representatives of $(N_G(T(F))\cap P_f)/P_T$; we deduce then immediately from the definitions (\cite[I.3.1.1 and I.8.1.1]{bt}) that $P_f$ is bounded; according to \cite[I.3.3.3]{bt}, there exists a maximal parahoric subgroup $P$ of $G$ containing $P_f$. Let $P_f^+$ be its pro-nilpotent radical; since $P$ contains $P_T$, there exists a concave function $f_P$ on $\Phi$ such that $P=P_{f_P}$; moreover, we have $P^+=U_{f'_P}$, where $f'_P$ is the concave function on $\Phi$ defined by $f'_P(\alpha)=1-f_P(-\alpha)$.
Since $P_f$ is contained in $P$, we have $f(\alpha)\geq f_P(\alpha)$ for every $\alpha$, which implies:
\[f(\alpha)\leq 1-f(-\alpha)\leq 1-f_P(-\alpha)=f'_P(\alpha)\]
for every $\alpha$; hence $P^+$ is contained in $P_f$. The group $P_f/P^+$ is therefore a closed subgroup of the reductive group $P/P^+$ containing the maximal torus $P_T/P^+$; moreover, if $\Phi_P$ is the root system of $P/P^+$ relatively to $P_T/P^+$, viewed as the subsystem of elements $\alpha\in\Phi$ such that $f_P(\alpha)+f_P(-\alpha)=0$, for every $\alpha\in\Phi_P$, we have $f'_P(\alpha)+f'_P(-\alpha)=2$, hence either $\alpha$ or $-\alpha$ is a root of $P_f/P^+$ relatively to $P_T/P^+$; moreover, we deduce easily from the commutator relations that the set of roots of $P_f/P^+$ is closed in $\Phi_P$. This set is then a parabolic subset of $\Phi_P$, and \cite[1.7. prop. 20]{bou} shows that $P_f/P^+$ is a parabolic subgroup of $P/P^+$; hence $P_f$ is a parahoric subgroup of $G(F)$, as required.

Assume now $f(\alpha)+f(-\alpha)=1$ for every $\alpha\in\Phi$. Then $P_f$ is pro-solvable, and $P_f/P^+$ is solvable, hence is a Borel subgroup of $P/P^+$; $P_f$ is then an Iwahori subgroup of $G(F)$, which shows the proposition. $\Box$

\section{Proof of the proposition \ref{nsgf}}

Assume now the conditions of the proposition \ref{nsgf} are satisfied. Let $H$ be an open bounded pro-nilpotent subgroup of $G(F)$ normalized by $P_T$; we'll show the existence of $T'$ and $f_H$.
 
We'll first show some preliminary results.

\begin{lemme}\label{chval}
It os possible to choose the valuation on $G(F)$ in such a way that $H$ is contained in $P_TG(F)_0^+$.
\end{lemme}

SInce $H$ is bounded, it is contained in some parahoric subgroup of $G$, hence the set $E_H$ of points of $\mathcal{B}$ fixed by $H$ is nonempty. Moreover, this set is bounded since $H$ is open, and stable by the action of $P_T$; hence by \cite[I.3.2.4]{bt}, there exists $x\in E_H$ such that $x$ is fixed by $P_T$, which is tantamount to say that the parahoric subgroup $G_x$ of $G(F)$ contains both $H$ and $P_T$.

Moreover, since $H$ is pro-solvable, the subgroup $H/H\cap G_{x,0}^+$ of the reductive group $G_x/G_{x,0}^+$ is a nilpotent subgroup normalized by the maximal torus $P_T/(P_T\cap G_{x,0}^+)$; it is then contained in the nilpotent radical of some Borel subgroup of $G_x/G_{x,0}^+$, hence $H$ is contained in the pro-nilpotent radical $I^+$ of the corresponding Iwahori subgroup $I$ of $G(F)$.

We may then assume $H=I^+$. We have $I=P_TU_{f_I}$, with $f_I$ being a concave function on $\Phi$ such that $f_I(\alpha)+f_I(-\alpha)=1$ for every $\alpha\in\Phi$. Let $P_0$ be a special parahoric subgroup containing $I$ and $P_1$ be its pro-nilpotent radical, let $\Delta_0$ be the basis of $\Phi$ associated to the Borel subgroup $I/P_1$ of $P_0/P_1$; assume the valuation on $G(F)$ has been chosen in such a way that we have $f_I(\alpha)=\frac 1{h(\Phi)}$ for every $\alpha$. Let $\beta$ be any element of $\Phi$; if $\beta$ is positive relatively to $\Delta_0$, we have $U_\beta(F)\cap I=U_\beta(F)\cap P_0$, and we obtain that $f_I(\beta)=\frac{h(\beta)}{h(\Phi)}>0$ ($h(\beta)$ being taken relatively to $\Delta_0$). If now $\beta<0$, we have:
\[f_I(\beta)=1-f_I(-\beta)=\frac{h(\Phi)-h(\beta)}{h(\Phi)}>0.\]
Since for every $\beta\in\Phi$, $f_I(\beta)>0$, $H$ is contained in $P_TG(F)_0^ +$, as required. $\Box$

For the sake of simplicity of notations, until the end of the proof of the proposition, we will assume $v_G(g)\in{\mth{Z}}\cup\{+\infty\}$ for every $g\in G(F)_0$; this can be done for example by multiplying the valuation previously considered by $h(\Phi$).

\begin{propo}\label{seprac2}
Assume at least one of the following conditions is fulfilled:
\begin{itemize}
\item $K$ has at least $7$ elements;
\item $K$ has $5$ elements and $\Phi$ has no component of type $A_1$, $C_n$  or $BC_n$;
\item $K$ has $4$ elements and $\Phi$ has no component of type $A_2$.
\end{itemize}
Let $h$ be any element of $H$; write:
\[h=\prod_{\alpha\in\overline{\Phi}}h_\alpha,\]
with $h_\alpha\in U'_\alpha$ for every $\alpha\in\overline{\Phi}$. Let $v=v_G(h)$; for every $\alpha\neq 0$ such that $v_\alpha(h_\alpha)=v$, $H\cap h_\alpha G_v^+$ is nonempty.
\end{propo}

Let $H_v$ (resp. $H_v^+$) be the subgroup of elements $h\in H$ such that $v(h)\geq v$ (resp. $v(h)>v)$; $H_v^+$ is a normal subgroup of $H_v$, and since $v>0$, the group $H_v/H_v^+$ is abelian, and is a subgroup of $G(F)_v/G(F)_v^+$, which is isomorphic to a direct product of copies of $K$; moreover, $H_v$ and $H_v^+$ are normalized by $P_T$, and if $P_T^+$ is the pro-nilpotent radical of $P_T$, we deduce from the commutator relations that we have:
\[[P_T^+,H_v]\subset H_v^+.\]
The torus $T(K)$ acts then on the abelian group $H_v/H_v^+$.
By a well-known result (see for example \cite[I.2.11(3)]{jan}), $H_v/H_v^+$ is the direct sum of its $T(K)$-weight subgroups; therefore, the assertion of the lemma is true as soon as all the elements of $\overline{\Phi}$ occuring as weights of $H_v/H_v^+$ are actually distinct characters of $T(K)$.

Let $\alpha\neq\beta$ be two elements of $\overline{\Phi}$; we'll show that $\alpha-\beta$ is nontrivial on $T(K)$. This is always the case if $K$ is infinite; if $K$ is finite and $q=\operatorname{card}(K)$, this is true as soon as $\alpha-\beta\not\in(q-1)X^*(S)$ (if and only if $\alpha-\beta\not\in(q-1)X^*(S)$ when $G$ is split).

Suppose first $\beta=0$, and let $\alpha^\vee$ be the coroot associated to $\alpha$; since by definition $\langle \alpha,\alpha^\vee\rangle =2$, $\alpha$ is nontrivial on $T(K)$ as soon as $q\geq 4$.

Suppose now $\beta=-\alpha$; by the same argument, $\alpha-\beta=2\alpha$ is nontrivial on $T(K)$ as soon as either $q=4$ or $q\geq 7$. Moreover, the case $\alpha\in 2X^*(T)$  occurs only when $\Phi$ has components of type $A_1$, $C_n$ or $BC_n$; in all other cases, there exists $\beta\in\Phi$ such that $\langle \alpha,\beta^\vee\rangle =1$, and $2\alpha$ is nontrivial on $T(K)$ as soon as $q\geq 4$.

Suppose now $\beta=2\alpha$; we have then $(2\alpha)^\vee=\frac 12\alpha^\vee$, and $\langle \alpha-\beta,(2\alpha)^\vee\rangle =1$; $\alpha-\beta=-\alpha$ is then nontrivial on $T(K)$ as soon as $q\geq 3$. By the same reasoning, if $\beta=-2\alpha$, $\alpha-\beta=3\alpha$ is nontrivial on $T(K)$ as soon as $q\geq 5$.

Suppose next $\alpha$ and $\beta$ are orthogonal. Since $\langle \alpha-\beta,\alpha^\vee\rangle =2$, $\alpha-\beta$ is nontrivial on $T(K)$ as soon as $q\geq 4$.

Suppose finally $\alpha$ and $\beta$ are non-orthogonal and linearly independant. According to \cite[7.5.1]{spr}, we then have either $\langle \alpha,\beta^\vee\rangle =\pm 1$ or $\langle \beta,\alpha^\vee\rangle =\pm 1$. Assume for example $\langle \beta,\alpha^\vee\rangle =\pm 1$; we obtain $0<\langle \alpha-\beta,\alpha^\vee\rangle \leq 3$, and $\alpha-\beta$ is nontrivial on $T(K)$ as soon as $q\geq 5$; moreover, if $\Phi$ is simply-laced and $\alpha$ and $\beta$ don't belong to any component of type $A_2$ of $\Phi$, there exists $\gamma\in\Phi$ such that $\langle \alpha-\beta,\gamma^\vee\rangle =1$, and we can replace $q\geq 5$ by $q\geq 4$. (Remember that the case $\Phi$ not simply-laced and $q=4$ is excluded by (G1).)

By summarizing the different cases, we get the assertion of the proposition. $\Box$

\begin{lemme}\label{sgrc}
Assume $p\neq 2$ or $F$ is absolutely unramified. Let $\alpha$ be an element of $\Phi$, let $v<v'$ be two integers such that there exists $h\in H$ such that $h=h_\alpha h_r$, with $h_\alpha\in U_\alpha(F)$, $v_\alpha(h_\alpha)=v$ and $v(h_r)\geq v'$. Then for any $h'_\alpha\in U_{\alpha,v}$, $H\cap h'_\alpha G_{v'}$ is nonempty.
\end{lemme}

Assume first $2\alpha\not\in\Phi$.
Set $h_\alpha=u_\alpha(x)$, $x\in F_\alpha$; for each $t\in P_T$, we have $th_rt^{-1}\in G_{v'}$, and $u_\alpha(\alpha(t)x)(th_rt^{-1})\in H$. If there exists $y\in\mathcal{O}_\alpha^*$ such that $t=\alpha^\vee(y)$, we obtain, setting $h'_r=th_rt^{-1}\in G_{v'}$:
\[u_\alpha(y^2x)h'_r\in H\]
for each $y$. Let $X$ be the subgroup of the elements $b$ of $\mathcal{O}_\alpha$ such that $u_\alpha(b)G_{v'}$ meets $H$; in order to show that $X=x\mathcal{O}_\alpha$, we only have to check that the ring $\mathcal{O}_\alpha$ is generated by the squares it contains.

Suppose first the characteristic $p$ of $K$ is odd or zero; we have:
\[z=\frac 12((z+1)^2-z^2-1),\]
which proves the assertion since $\frac 12=2(\frac 1{2^2})$.

Suppose now $p=2$. Since $K_\alpha$ is perfect, any element of $K_\alpha^*$ is a square; moreover, we have $2=1^2+1^2$. We then only have to show that every element of $1+{\mathfrak{p}}_\alpha$ belongs to the subring of $\mathcal{O}_\alpha$ generated by the squares, which is simply done by remarking that if $z\in 1+{\mathfrak{p}}_\alpha$, since $F_\alpha$ is absolutely unramified, $\frac{z+1}2$ and $\frac{z-1}2$ belong to $\mathcal{O}_\alpha$ and $z=(\frac{z+1}2)^2-(\frac{z-1}2)^2$.

(Remark: when $p=2$, the result is false when $F$ is absolutely ramified. It is easy to check for example that for $F=F_0[\sqrt{2}]$, where $F_0$ is any unramified extension of ${\mth{Q}}_2$, and $G=SL_2$, we can obtain counterexamples to the assertion of this lemma, and even to the proposition \ref{nsgf}.)

Assume now $2\alpha\in\Phi$; since $\Phi$ is then not simply-laced, by (G1) we must have $p\neq 2$. By an explicit computation in $SU_3$, we easily see that for an appropriate choice of the isomorphism between $H_\alpha$ and $U_\alpha(F)$, we have, setting again $t=\alpha^\vee(y)$:
\[tu_\alpha(x)t^{-1}=u_\alpha(y^2\sigma(y)^{-1}x),\]
with $\sigma$ being the nontrivial element of $\operatorname{Gal}(F_\alpha/F_{2\alpha})$. We claim that the elements $y^2\sigma(y)^{-1}$, $y\in\mathcal{O}_\alpha$; generate the ring $\mathcal{O}_\alpha$: for $y\in\mathcal{O}_{2\alpha}$, we simply have $y^2\sigma(y)^{-1}=y$, and for any $y\in\mathcal{O}_\alpha^*$ such that $\sigma(y)=-y$ (such an $y$ exists because $F_\alpha/F_{2\alpha}$ is unramified and $p\neq 2$), we have $y^2\sigma(y)^{-1}=-y$. Since such an element and $1$ generate $\mathcal{O}_\alpha$ as a $\mathcal{O}_{2\alpha}$-module, the claim s proved.

On the other hand, we have:
\[[u_\alpha(x')U_{2\alpha,2c}G_{v'},u_\alpha(x'')U_{2\alpha,2c}G_{v'}]\subset u_{2\alpha}(\sigma(x')x''-x'\sigma(x''))G_{v'},\]
where $\sigma$ is the nontrivial element of $\operatorname{Gal}(F_\alpha/F_{2\alpha})$; moreover, if $x',x''$ are chosen such that the image of $\frac{x'}{x''}$ in $K_\alpha$ doesn't belong to $K$, the valuation of $\sigma(x')x''-x'\sigma(x'')$ is exactly $2c$;
since $2(2\alpha)\not\in\Phi$, we deduce from the preceding case that for each $y\in\varpi^{2c}\mathcal{O}_{2\alpha}$, $H$ meets $u_{2\alpha}(y)G_{v'}$. We combine these two assertions to obtain that $H$ meets $h_\alpha G_{v'}$ for every $h_\alpha\in U_\alpha(F)$ such that $v_\alpha(h_\alpha)\geq c$, as required. $\Box$

For every $\alpha\in\Phi$ and every positive elements $v\leq c$ of ${\mth{Z}}$, let $f_{\alpha,c,v}$ be the function on $\Phi$ such that $f_{\alpha,c,v}(\alpha)=v$, $f_{\alpha,c,v}(2\alpha)=v$ if $2\alpha\in\Phi$, and $f_{\alpha,c,v}(\beta)=c$ for every $\beta\neq\alpha,2\alpha$.

\begin{lemme}\label{facv}
The function $f_{\alpha,c,v}$ is a concave function; moreover, $U_{f_{\alpha,c+1,v+1}}$ is normal in $U_{f_{\alpha,c,v}}$ and the quotient $U_{f_{\alpha,c,v}}/U_{f_{\alpha,c+1,v+1}}$ is abelian.
\end{lemme}

Let $\beta,\gamma$ be two elements of $\phi$ such that $\beta+\gamma\in\phi$. If $\beta$ and $\gamma$ are both contained in $\{\alpha,2\alpha\}$, then $\beta=\gamma=\alpha$ and $\beta+\gamma=2\alpha$, and $f_{\alpha,c,v}(\beta+\gamma)=v\leq 2v=f_{\alpha,c,v}(\beta)+f_{\alpha,c,v}(\gamma)$; if now at least one of them is different from $\alpha$ and $2\alpha$, we have:
\[f_{\alpha,c,v}(\beta)+f_{\alpha,c,v}(\gamma)\geq c+v\geq c\geq f_{\alpha,c,v}(\beta+\gamma).\]
since it is obvious that $f_{\alpha,c,v}(\beta)+f_{\alpha,c,v}(-\beta)>0$ for each $\beta\in\phi$, the concavity of $f_{\alpha,c,v}$ is proved; moreover, the above inequalities imply:
\[[U_{\alpha,f_{\alpha,c,v}(\alpha)},U_{\alpha,f_{\alpha,c,v}(\alpha)}]\subset U_{2\alpha,2v}\subset U_{2\alpha,f_{\alpha,c+1,v+1}(2\alpha)}\]
if $2\alpha\in\phi$, and:
\[[U_{\beta,f_{\alpha,c,v}(\beta)},U_{\gamma,f_{\alpha,c,v}(\gamma)}]\subset U_{\beta+\gamma,c+v}\subset U_{\beta+\gamma,f_{\alpha,c+1,v+1}\beta+\gamma)}\]
for every $\beta,\gamma\in\phi$ such that $\beta+\gamma\in\phi$ and either $\beta$ or $\gamma$ doesn't belong to $\{\alpha,2\alpha\}$; moreover, we have, for every $\beta\in\phi$:
\[[U_{\beta,f_{\alpha,c,v}(\beta)},U_{0,c}]\subset U_{\beta,c+v}\subset U_{\beta,f_{\alpha,c+1,v+1}(\beta)}\]
and, since either $\beta$ or $-\beta$ doesn't belong to $\{\alpha,2\alpha\}$:
\[[U_{\beta,f_{\alpha,c,v}(\beta)},U_{-\beta,f_{\alpha,c,v}(-\beta)}]\subset U_{0,c+v}\subset U_{0,c+1};\]
We deduce from all these inclusions that $U_{f_{\alpha,c+1,v+1}}$ is normal in $U_{f_{\alpha,c,v}}$ and that the quotient $U_{f_{\alpha,c,v}}/U_{f_{\alpha,c+1,v+1}}$ is abelian, as required. $\Box$

Consider now the group $G(F)_{c,v}=U_{0,v}G(F)_c$; we have the following result:

\begin{lemme}\label{facv2}
The group $G(F)_{c+1,v+1}$ is normal in $G(F)_{c,v}$, and the quotient $G(F)_{c,v}/G(F)_{c+1,v+1}$ is abelian.
\end{lemme}

The proof is analoguous to the proof of the previous lemma. $\Box$

Now we'll prove proposition \ref{nsgf}. According to the remarks made in the previous section and to the lemma \ref{sgrc}, we only have to show that we have:
\[H=\prod_{\alpha\in\overline{\Phi}}(H\cap U_\alpha(F)),\]
the product being taken in any order.

Let $h=\prod_{\alpha\in\overline{\Phi}}h_\alpha$ be defined as in the proposition \ref{seprac2}, and set $v=v_G(h)$; assume moreover that the product is chosen according to some order on $\overline{\Phi}$ satisfying the following conditions:
\begin{itemize}
\item if $\alpha$ is any element of $\Phi$, $0\leq\alpha$;
\item if $\alpha,\beta$ are elements of $\Phi$, we have $\alpha<\beta$ if and only if $v(h_\alpha)>v(h_\beta)$.
\end{itemize}
With the help of the commutator relations, it is easy to check that this is always possible.

Since $H$ is open, there exists an integer $c_0$ such that $G_{c_0}\subset H$. We'll claim that for every $\alpha\in\overline{\Phi}$, $H$ contains $h_\alpha$. This may be done by proving that $H$ contains elements of $h_\alpha G(F)_c$, with $c$ arbitrarily large; since this will in particular be true for $c=c_0$, we will obtain that the element $h_\alpha$ itself belongs to $H$. We'll show the claim by induction on $c_0-v$, the case $v=c_0$ being trivial.

First assume $\alpha\in\Phi$ and $v_\alpha(h_\alpha)=v$; we'll proceed by induction on $c$. The case $c=v+1$ is simply the result of the proposition \ref{seprac2}; assume now $c>v+1$, and let $h_{c-1}$ be an element of $h_\alpha G(F)_{c-1}\cap H$; such an element exists by the induction hypothesis. The element $h_{c-1}$ belongs to $U_{f_{\alpha,c-1,v}}$; we deduce then from the lemma \ref{facv}, by using as in the proof of proposition \ref{seprac2} the fact that the quotient $U_{f_{\alpha,c-1,v}}/U_{f_{\alpha,c,v+1}}$ can be viewed as a rational representation of $T(K)$, that $h_\alpha U_{f_{c,v+1}}\cap H$ is nonempty. Let $h_c$ be any element of this intersection; the lemma \ref{sgrc} allows us to choose $h_c$ in $h_\alpha G(F)_c$, as required.

Assume $v_0(h_0)>v$. By an easy induction on $v'$, we see that the element:
\[h'=h(\prod_{\alpha,v_\alpha(h_\alpha)=v}h_\alpha)^{-1}=\prod_{\beta\in\overline{\Phi},v_\beta(h_\beta)>v}h_\beta\]
belongs to $H$; our claim follows then from the induction hypothesis applied to $h'$.

Now assume $v=v_0(h_0)$. Since $H$ then meets $G(F)_{v,c_0}$, we deduce from the lemma \ref{facv2} and the proof of the proposition \ref{seprac2}, as above, that $H$ contains an element of $h'_0G(F)_{c_0}$, with $h'_0$ being an element of $T(F)$ such that $h'_0{}^{-1}h\in U_{0,v+1}$; hence $h'_0$ belongs to $H$. Since $h'_0{}^{-1}h\in H$ and $v_G(h'_0{}^{-1}h)>v$, we can use the induction hypothesis again.

Since the above claim is true for any $h\in H$, we have just shown the inclusion:
\[H\subset\prod_{\alpha\in\overline{\Phi}}(H\cap U_\alpha(F));\]
the other inclusion being obvious, the proposition is proved. $\Box$

\section{Normal functions}

From now on and until the end of the paper, $H$ will be a normal open pro-nilpotent subgroup of some parahoric subgroup of $G(F)$. Moreover, until the end of the proof of the theorem \ref{cln}, $G$ will be assumed to be quasi-simple, which amounts to say that $\Phi$ is connected.

Since every parahoric subgroup of $G(F)$ contains an Iwahori subgroup, and all Iwahori subgroups of $G(F)$ are conjugated, we may without loss of generality choose one of them and only consider its normal open pro-nilpotent subgroups; we will then make the following assumptions:
\begin{itemize}
\item the (extended) valuation $(v_\alpha)_{\alpha\in{\overline{\Phi}}}$ has been chosen in such a way that the subgroup $G(F)_0$ of $G(F)$ is a special parahoric subgroup. Moreover, for every $\alpha\in\Phi$, $v_\alpha(U'_\alpha)={\mth{Z}}\cup\{+\infty\}$;
\item let $B(K)$ be the Borel subgroup of $G(K)\simeq G(F)_0/G(F)_0^+$ associated to $\Delta$; $I$ is the inverse image of $B(K)$ in $G(F)_0$.
\end{itemize}

We then have $I=P_TU_{f_I}$, where $f_I$ is the fonction on $\Phi$ defined by $f_I(\alpha)=0$ (resp. $f_I(\alpha)=1$) when $\alpha$ is a positive (resp. negative) root relatively to $\Delta$, or equivalently $f_I(\alpha)=\varepsilon_\alpha$ for every $\alpha\in\Phi$.

We may easily check that if $\alpha,\beta$ are elements of $\Phi$ such that $\alpha+\beta\in\Phi$, we have:
\begin{itemize}
\item if $\alpha+\beta$ is greater than $\alpha$ and $\beta$, then $f_I(\alpha+\beta)=f_I(\alpha)+f_I(\beta)$;
\item if $\alpha+\beta$ is lesser than $\alpha$ and $\beta$, then $f_I(\alpha+\beta)=f_I(\alpha)+f_I(\beta)-1$.
\end{itemize}

Write $H=T'U_{f_H}$ as in the preceding proposition. According to \cite[I.6.4.43]{bt}, $f_H$ must satisfy:
\[f_H(\alpha+\beta)\leq f_I(\alpha)+f_H(\beta)\]
for every $\alpha,\beta\in\Phi$ such that $\alpha+\beta\in\Phi$. Conversely, for any $f$ satisfying the above condition, there exists $T'\subset T$ such that $T'U_f$ is normal in $I$. We'll say $f$ is a normal function for $I$ if it satisfies that condition.

Consider the application $f_\varepsilon=f-f_I$ from $\Phi$ to ${\mth{Z}}$; we deduce from the above inequality that for every $\alpha,\alpha'\in\Phi$ such that $\alpha+\alpha'\in\Phi$ and $\alpha+\alpha'\geq\alpha$, we have:
\[f_\varepsilon(\alpha+\alpha')\leq f_\varepsilon(\alpha).\]
With the proposition \ref{decprg} and an easy induction, we obtain that for every $\alpha,\beta\in\Phi$ such that $\alpha\geq\beta$, $f_\varepsilon(\alpha)\leq f_\varepsilon(\beta)$. Moreover, we have the following results:

\begin{lemme}
For every $\alpha\in\Phi$, we have $f_\varepsilon(\alpha)\leq f_\varepsilon(-\alpha)+1$.
\end{lemme}

Consider the element $u_\alpha(\varpi^{f_I(\alpha)})$ of $I$; assuming $2\alpha\not\in\Phi$, we have:
\[u_\alpha(\varpi^{f_I(\alpha)})u_{-\alpha}(\varpi^{f_H(-\alpha)})u_\alpha(-\varpi^{f_I(\alpha)})\]
\[=u_{-\alpha}(\varpi^{f_H(-\alpha)})\alpha^\vee(1+\varpi^{f_H(-\alpha)+f_I(\alpha)}x)u_\alpha(\varpi^{(f_H(-\alpha)+2f_I(\alpha)}x'),\]
\[=u_{-\alpha}(\varpi^{f_H(-\alpha)})\alpha^\vee(1+\varpi^{1+(f_\varepsilon)(-\alpha)}x)u_\alpha(\varpi^{1+(f_\varepsilon)(-\alpha)+f_I(\alpha)}x'),\]
where $x$ and $x'$ are elements of $1+{\mathfrak{p}}$. We deduce then from the proposition \ref{nsgf} that $U_{\alpha,(1+(f_\varepsilon)(-\alpha)+f_I(\alpha)}\subset H$, which implies $f_H(\alpha)\leq 1+(f_\varepsilon)(-\alpha)+f_I(\alpha)$; hence the result. The proof when $2\alpha\in\Phi$ is left to the reader. $\Box$

\begin{lemme}
Let $\alpha$ be an element of $\Phi$; set $v=\operatorname{Sup}(f_\varepsilon(\alpha),f_\varepsilon(-\alpha))$. For every $\beta\in\Phi$, $f_\varepsilon(\beta)\leq v+1$.
\end{lemme}

Assume there exists $\beta$ such that $f_\varepsilon(\beta)\geq v+2$; since we can always replace $\beta$ by a smaller element of $\Phi$, we may assume $\beta\in\Delta'$. Since $(\alpha,\varepsilon_\alpha)+(-\alpha,\varepsilon_{-\alpha})=(0,1)\geq(\beta,\varepsilon_\beta)$, we have either $\beta\leq\alpha$ or $\beta\leq-\alpha$; we may assume for example $\beta\leq\alpha$, which implies $-\alpha\leq-\beta$ and $f_\varepsilon(-\beta)\geq f_\varepsilon(\beta)-1\geq v+1>f_\varepsilon(-\alpha)$, hence a contradiction. $\Box$

By a similar reasoning, we obtain:

\begin{lemme}
Let $\alpha$ be any element of $\Phi$; set $v=\operatorname{Inf}(f_\varepsilon(\alpha),f_\varepsilon(\beta))$ Then for every $\beta\in\Phi$, $f_\varepsilon(\beta)\geq v-1$.
\end{lemme}

We deduce from these three lemmas the following corollary:

\begin{coro}
There exists $v'\in{\mth{N}}$ such that $\operatorname{Im}(f_\varepsilon)$ is contained in $\{v',v'+1,v'+2\}$. Moreover, if $\alpha,\beta\in\Phi$ are such that $f_\varepsilon(\alpha)=v'$ and $f_\varepsilon(\beta)=v'+2$, then $f_\varepsilon(-\alpha)=f_\varepsilon(-\beta)=v'+1$.
\end{coro}

Set $v'=\operatorname{Inf}_{\alpha\in\Phi}f_\varepsilon(\alpha)$; let $\alpha$ be an element of $\Phi$ such that $f_\varepsilon(\alpha)=v'$. Then according to the first lemma, $f_\varepsilon(-\alpha)\leq v'+1$, and according to the second one, $f_\varepsilon(\beta)\leq v'+2$ for every $\beta\in\Phi$. The second assertion of the corollary is an immediate consequence of the second and third lemmas. $\Box$

In the sequel, we'll say $f_H$ is:
\begin{itemize}
\item of type $1$ if $f_\varepsilon$ either is constant or takes onky two consecutive different values $v',v'+1$;
\item of type $2$ else.
\end{itemize}

\section{$\Delta'$-complete subsets of $\Phi$}

Consider the subset $\Psi$ of the elements $\alpha\in\Phi$ such that $f_\varepsilon(\alpha)$ is maximal; we deduce from the previous section that for every $\alpha\in\Psi$ and every $\beta\in\Phi$ such that $\beta\leq\alpha$, $\beta\in\Psi$. We'll say a subset of $\Phi$ satisfying that condition is a $\Delta'$-complete subset of $\Phi$. This section will be devoted to the study of some properties of such subsets, which will be useful in the determination of the $x$ and $r$ of the theorem.

Assume for example $f_H$ is of type $1$; $x$ and $r$ must satisfy $\alpha(x)-r\leq-v'-\varepsilon_\alpha$ for every $\alpha\in\Psi$, and $\alpha(x)-r\geq-v'-\varepsilon_\alpha$ for every $\alpha\not\in\Psi$. Hence if $\alpha_1,\dots,\alpha_t$ (resp. $\beta_1,\dots,\beta_s$) are elements of $\Psi$ (resp. $\Phi-\Psi$) whose sum is zero, we must have, setting $z=r-v'$:
\[\sum_{i=1}^t\varepsilon_{\alpha_i}\leq tz;\]
\[\sum_{j=1}^s\varepsilon_{\beta_j}\geq sz.\]
The purpose of the proposition \ref{smph} is to prove the existence of some $z$ satisfying the above conditions; when this is done, the proposition \ref{ccz} ensures the existence of a suitable $x$.

Until the end of the section, $\Psi$ will be any $\Delta'$-complete subset of $\Phi$. We have:

\begin{propo}\label{smph}
Let $\alpha_1,\dots,\alpha_t$ (resp. $\beta_1,\dots,\beta_s$) be elements of $\Psi$ (resp. $\Phi-\Psi$) such that $\sum_{i=1}^t(\alpha_i,\varepsilon_{\alpha_i})=(0,c)$ and $\sum_{j=1}^s(\beta_j,\varepsilon_{\beta_j})=(0,d)$ for some integers $c,d$.
Assume $s,t>0$; then $\frac ct\leq\frac ds$, and if $\frac ct=\frac ds$, for every $i,j$, $\alpha_i$ and $\beta_j$ are strongly orthogonal.
\end{propo}

As a first remark, we see that we can assume $\sum_{i=1}^{t'}\alpha_i\in\overline{\Phi}$ for every $t'\leq t$; let's proceed by induction on $t'$, the case $t'=1$ being obvious. If the sum is zero, then $\sum_{i=1}^{t'+1}\alpha_i=\alpha_{t'+1}\in\overline{\Phi}$, as required; if the sum is nonzero, since $\sum_{i=t'+1}^t\alpha_i=-\sum_{i=1}^{t'}\alpha_i$, there exists $i'>t'$ such that $(\sum_{i=1}^{t'+1}\alpha_i,\alpha_{i'})<0$, and by rearranging the $\alpha_i$, $i>t$, we may assume $i'=t+1$, hence the result. (This is a slight variant of \cite[I.proposition 1.19]{bou} but not a direct consequence of it.) From now on and unless another rearranging is explicitly mentioned, we will assume the $\alpha_i$ follow that property; we will make the same assumption about the $\beta_j$.

First we'll prove we can assume $\Phi$ to be simply-laced. (This shouldn't be really necessary to run the rest of the proof, but it makes it somewhat simpler.) It is well-known (see for example \cite[par. 11]{st}, and in particular theorem 32) that when $\Phi$ is not simply-laced, there exists a simply-laced root system $(\Phi_{sl},V_{sl})$ and an automorphism $\sigma$ of this root system, satisfying the following conditions:
\begin{itemize}
\item there exists a basis $\Delta_{sl}$ of $\Phi_{sl}$ such that $\sigma(\Delta_{sl})=\Delta_{sl}$;
\item $V$ can be identified to the subspace of the elements of $V_{sl}$ fixed by $\sigma$; 
\item let $W_{sl}$ be the Weyl group of $\Phi_{sl}$, let $(.,.)_{sl}$ be a $W_{sl}$-invariant scalar product, and let $\pi$ be the orthogonal projection, according to this scalar product, from $V_{sl}$ to $V$. Then $\pi(\Phi_{sl})=\Phi$ and $\pi(\Delta_{sl})=\Delta$; moreover, $W$ is the subgroup of the elements of $W_{sl}$ commuting with $\sigma$, and $(.,.)$ is the restriction of $(.,.)_{sl}$ to $V$.
\end{itemize}
Let $\alpha$ be an element of $\Phi$ and let $\alpha_{sl}$ be an element of $\Phi_{sl}$ whose image in $\Phi$ is $\alpha$; we have:
\[\alpha=\frac 1d\sum_{i=0}^{d-1}\sigma^i(\alpha_{sl}),\]
where $d$ is the order of $\sigma$. Moreover, the image of any positive element of $\Phi_{sl}$ is a positive element of $\Phi$, hence $\varepsilon_{\alpha_{sl}}=\varepsilon_{\pi(\alpha_{sl})}$ for every $\alpha_{sl}\in\Phi_{sl}$.

Let $\Phi_{sl,aff}$ be the affine root system associated do $\Phi_{sl}$; the morphism $\pi$ extends canonically to a morphism $\Phi_{sl,aff}\mapsto\Phi_{aff}$. Moreover, if $\Delta'_{sl}$ is the extended simple root system of $\Phi_{sl}$ associated do $\Delta_{sl}$, we have $\pi(\Delta'_{sl})=\Delta'$, and for every $(\alpha_{sl},v)\in\Phi^+_{sl,aff}$, $h(\pi(\alpha_{sl},v))=h(\alpha_{sl},v)$; in particular, $h(\Phi_{sl})=h(\Phi)$, and $(\alpha_{sl},v)\leq(\beta_{sl},v')$ if and only if $\pi(\alpha_{sl},v)\leq\pi(\beta_{sl},v')$. We finally have the following results:

\begin{lemme}
Let $\alpha,\beta$ be elements of $\Phi$ such that $\alpha+\beta\in\Phi$; there exist $\alpha_{sl},\beta_{sl}\in\Phi_{sl}$ such that $\alpha_{sl}+\beta_{sl}\in\Phi_{sl}$, $\pi(\alpha_{sl})=\alpha$ and $\pi(\beta_{sl})=\beta$.
\end{lemme}

Assume first $(\alpha,\beta)<0$, and let $\alpha'_{sl}$ (resp. $\beta'_{sl}$) be any element of $\Phi_{sl}$ whose image in $\Phi$ is $\alpha$ (resp. $\beta$). We have:
\[(\sum_{i=0}^{d-1}\sigma^i(\alpha'_{sl}),\sum_{j=0}^{d-1}\sigma^j(\beta'_{sl}))<0;\]
there exist then $i,j$ such that $(\sigma^i(\alpha'_{sl}),\sigma^j(\beta_{sl}))<0$, which proves the result. Suppose now $(\alpha,\beta)\geq 0$; we then have $(\alpha+\beta,-\beta)<0$, and the above reasoning applied to $\alpha+\beta$ and $-\beta$ yields elements $(\alpha+\beta)_{sl}$ and $(-\beta)_{sl}$ of $\Phi$ whose sum is an element $\alpha_{sl}$ such that $\pi(\alpha_{sl})=\alpha$, hence again the result. $\Box$

\begin{coro}
Let $\alpha_1,\dots,\alpha_t$ be elements of $\Phi$ such that $(\gamma,v)=\sum_{i=1}^t(\alpha_i,\varepsilon_{\alpha_i})$ is an element of $\overline{\Phi}_{aff}$. There exist $\alpha_{sl,1},\dots,\alpha_{sl,t},\gamma_{sl}\in\Phi$ such that $\pi(\alpha_{s_i})=\alpha$ for every $i$, $\pi(\gamma_{sl})=\gamma$ and $(\gamma_{sl},v)=\sum_{i=1}^t(\alpha_{sl,i},\varepsilon_{\alpha_{sl,i}})$.
\end{coro}

This corollary follows from the previous lemma and an easy induction. $\Box$

Assume the proposition is true in the simply-laced case. Let $\Psi_{sl}$ be the subset of the elements of $\Phi_{sl}$ whose image belongs to $\Psi$; since $\pi$ preserves the partial order, $\Psi_{sl}$ is a $\Delta'$-complete subset of $\Phi_{sl}$. Let $\alpha_{sl,1},\dots,\alpha_{sl,t}$ (resp. $\beta_{sl,1},\dots,\beta_{sl,s}$) be elements of $\Phi_{sl}$ whose images are the $\alpha_i$ (resp. the $\beta_j$) and which satisfy the conditions of the previous corollary; the fact that $\frac ct\leq\frac ds$ is simply the first assertion in the simply-laced case; moreover, if for some $i,j$, $\alpha_i$ and $\beta_j$ are not strongly orthogonal, by eventually conjugating all the $\beta_j$ by $\sigma^k$ for some $k$, we can assume $\alpha_{sl,i}$ and $\beta_{sl,i}$ are not strongly orthogonal, and by the simply-laced case, we obtain $\frac ct<\frac ds$, as required.

We will now assume $\Phi$ is simply-laced. In this case, two elements of $\Phi$ are strongly orthogonal if and only if they are orthogonal.

First we'll observe that if there exist families $(\alpha_1,\dots,\alpha_t)$ and $(\beta_1,\dots,\beta_s)$ such that $\frac ct\leq\frac ds$, then the families consisting of respectively $d$ copies of the $\alpha_i$ and $c$ copies of the $\beta_j$ satisfy $t\leq s$. Hence if there don't exist any families $(\alpha_1,\dots,\alpha_t),(\beta_1,\dots,\beta_s)$ such that $c=d$ and $t\leq s$, then for every families, we have $\frac ct>\frac ds$ and the proposition is proved. We'll then assume such families do actually exist.

Now we'll make the following claim: if $(\alpha_1,\dots,\alpha_t)$ and $(\beta_1,\dots,\beta_t)$ satisfy the above condition, and are such that $c$ is minimal among all families satisfying it, then $t=s$ and for every $i,j$, $\alpha_i$ and $\beta_j$ are orthogonal. First we'll show:

\begin{lemme}
For every proper nonempty subset $J$ of $\{1,\dots,s\}$, $\sum_{j\in J}\beta_j\neq 0$.
\end{lemme}

Assume there exists some $J$ such that $\sum_{j\in J}\beta_j=0$; let $s'$ be the cardinal of $J$, and set $c'=\sum_{j\in J}\varepsilon_{\beta_j}$; we have $0<c'<c$. Let $t'$ be an element of $\{1,\dots,t\}$ such that $\sum_{i=1}^{t'}(\alpha_i,\varepsilon_{\alpha_i})\geq(0,c')$ and $\sum_{i=1}^{t'-1}(\alpha_i,\varepsilon_{\alpha_i})<(0,c')$; according to the proposition \ref{infsum}, there exist $\alpha'_1,\dots,\alpha'_{t''-1}\in\Psi\cup\{0\}$, with $t''\leq t'$, such that:
\[\sum_{i=1}^{t''}(\alpha'_i,\varepsilon_{\alpha'_i})=(0,c');\]
by minimality of $c$, we must have $t''>s'$, hence $t'\geq s'+1$. Moreover, the same argument applied to $\alpha_{t'},\dots,\alpha_{t}$ yields $t-t'+1\geq s-s'+1$; we then obtain $t+1\geq s+2$, which contradicts our assumptions. Hence the lemma. $\Box$

Now we'll prove the claim. Suppose there exists $i,j$ such that $(\alpha_i,\beta_j)>0$; we may assume $i=t$ and $j=s$. We have $\beta_s>\alpha_t$, $\delta_1=\alpha_t-\beta_s\in\Phi$ and:
\[\gamma=\sum_{j=1}^{s-1}(\beta_j,\varepsilon_{\beta_j})=\sum_{i=1}^{t-1}(\alpha_i,\varepsilon_{\alpha_j})-(\delta_1,\varepsilon_{\delta_1}).\]
According to the previous lemma, $\gamma$ is nonzero; there exists then $j\in\{1,\dots,s-1\}$ such that $(\beta_j,\gamma)>0$, hence either $(\beta_j,\alpha_i)>0$ for some $i$ or $(\beta_j,-\delta_1)>0$. Assume $j=s-1$; in the first case, assuming $i=t-1$, we obtain;
\[\sum_{j=1}^{s-2}(\beta_j,\varepsilon_{\beta_j})=\sum_{i=1}^{t-2}(\alpha_i,\varepsilon_{\alpha_j})-(\delta_1,\varepsilon_{\delta_1})-(\delta_2,\varepsilon_{\delta_2}),\]
with $\delta_2=\beta_{s-1}-\alpha_{t-1}$, and in the second case:
\[\sum_{j=1}^{s-2}(\beta_j,\varepsilon_{\beta_j})=\sum_{i=1}^{t-1}(\alpha_i,\varepsilon_{\alpha_j})-(\delta_1,\varepsilon_{\delta_1})-(\beta_{s-1},\varepsilon_{\beta_{s-1}}).\]
By iterating the process, we finally obtain, after eventually rearranging the $\beta_j$:
\[0=\sum_{i=1}^{t'}(\alpha_i,\varepsilon_{\alpha_j})-\sum_{k=1}^u(\delta_k,\varepsilon_{\delta_k})-\sum_{j=1}^{s'}(\beta_j,\varepsilon_{\beta_j}),\]
with $s'+u=s$ and $t'+u=t$, hence $t'\leq s'$. Consider the equality:
\[\sum_{k=1}^u(\delta_k,\varepsilon_{\delta_k})=\sum_{i=1}^{t'}(\alpha_i,\varepsilon_{\alpha_i})+\sum_{j=1}^{s'}(-\beta_j,-\varepsilon_{\beta_j}).\]
Although the members of this equality are not necessarily affine roots, we can here use a similar reasoning as in the proof of the proposition \ref{infsum}: if the sum $\delta$ of the $\delta_k$ is nonzero, there exists an element $\gamma\in\Delta'$ such that $(\delta,\gamma)>0$, hence $(\delta_i,\gamma)>0$ for some $i$; by subtracting $(\gamma,\varepsilon_\gamma)$ to some appropriate term of the right-hand side and iterating, we finally obtain an equality such as:
\[(0,c')=\sum_{i=1}^{t''}(\alpha'_i,\varepsilon_{\alpha'_i})-\sum_{j=1}^{s''}(\beta'_j,\varepsilon_{\beta'_j})-(0,s'-s''),\]
with $t'\leq t$, $s''\leq s$ and the $\alpha'_i$ (resp. the $\beta'_j$) being elements of $\Psi$ (resp. $\Phi-\Psi$); the last term of the right-hand side corresponds to the $(-\beta_j,-\varepsilon_{\beta_j})$ which are reduced to $(0,-1)$ that way. Hence:
\[\sum_{i=1}^{t''}(\alpha'_i,\varepsilon_{\alpha'_i})=\sum_{j=1}^{s''}(\beta'_j,\varepsilon_{\beta'_j})+(0,c'+s'-s'').\]
Moreover, we have $t''\leq t'\leq s'\leq s''+(c'+s'-s'')$. If $(\alpha'_i,\beta'_j)>0$ for some $i,j$, we can iterate the whole process to get even smaller sums satisfying similar inequalities.

Assume now $(\alpha'_i,\beta'_j)=0$ for every $i,j$; we will show by induction on $c'+s'-s''$ that it leads to a contradiction. First remark that the sum of the $\alpha'_i$ and the sum of the $\beta'_j$ must be zero, since they are equal to each other and their product is zero. If $c'+s'-s''=0$, since $\sum_{i=1}^{t''}(\alpha'_i,\varepsilon_{\alpha'_i})<(0,c)$, the equality is impossible by minimality of $c$; if now $c'+s'-s''>0$, by subtracting $(\alpha'_{t''},\varepsilon_{\alpha'_{t''}})$ to both sides, we obtain;
\[\sum_{i=1}^{t''-1}(\alpha'_i,\varepsilon_{\alpha'_i})=\sum_{j=1}^{s''}(\beta'_j,\varepsilon_{\beta'_j})+(-\alpha'_{t''},c'+s'-s''-1+\varepsilon_{-\alpha'_{t''}}).\]
after replacing $(-\alpha'_{t''},c'+s'-s''-1+\varepsilon_{-\alpha'_{t''}})$ by $(0,c'+s'-s''-1)$ in the right-hand side and applying the proposition \ref{infsum} to the left-hand side, we can use the induction hypothesis to obtain the desired contradiction.

Since assuming $(\alpha_i,\beta_j)>0$ for some $i,j$ leads to a contradiction, and since the sum of the $\alpha_i$ is zero, we have just shown $(\alpha_i,\beta_j)=0$ for every $i,j$. Let now $\delta$ be an element of $\Delta'$ such that $(\alpha_t,\delta)>0$; since the sum of the $(\beta_j,\varepsilon_{\beta_j})$ is greater than $(\delta,\varepsilon_\delta)$, there exists $j$ such that $\delta\leq\beta_j$. Assuming $j=s$, we have:
\[\sum_{j=1}^{s-1}(\beta_j,\varepsilon_{\beta_j})<\sum_{i=1}^{t-1}(\alpha_i,\varepsilon_{\alpha_i})+(\alpha_t-\delta),\]
hence, by applying the proposition \ref{infsum}:
\[\sum_{j=1}^{s-1}(\beta_j,\varepsilon_{\beta_j})=\sum_{i=1}^t(\alpha'_i,\varepsilon_{\alpha'_i}),\]
the $\alpha'_i$ being elements of $\Psi\cup\{0\}$. If $t<s$, a similar reasoning as above leads to a similar contradiction; hence $t=s$ and the claim is proved.

Let's write $c_0$ for the minimal $c$ defined as before, and set $t_0=\frac{c_0}{z(\Psi)}$; for convenience, we will also set $s_0=t_0$ and $d_0=c_0$. We have the following result:

\begin{lemme}\label{relp}
The integers $c_0$ and $t_0$ are relatively prime.
\end{lemme}

Assume they are not, and let $c',t'$ be positive and relatively prime integers such that $\frac{c'}{t'}=\frac{c_0}{t_0}$. We have either $\sum_{i=1}^{t'}(\alpha_i,\varepsilon_{\alpha_i})\geq(0,c')$ or $\sum_{i=t'+1}^{t_0}(\alpha_i,\varepsilon_{\alpha_i})\geq(0,c_0-c')$. Moreover, $c_0$ is a multiple of $c'$; by replacing, in the first case, the family $(\alpha_1,\dots,\alpha_{t'})$ by a family made of $\frac{c_0-c'}{c'}$ copies of it, we are reduced to the second case. By applying proposition \ref{infsum}, we then obtain $\sum_{i=t'+1}^{t_0}(\alpha'_i,\varepsilon_{\alpha'_i})=(0,c_0-c')$, the $\alpha'_i$ being elements of $\Psi\cup\{0\}$.

Similarly (setting $s'=t'$ for convenience), either $\sum_{j=1}^{s'}(\beta_j,\varepsilon_{\beta_j})\leq(0,c')$ or $\sum_{j=s'+1}^{s_0}(\beta_j,\varepsilon_{\beta_j})\leq(0,c_0-c')$. By the same argument as above, we only have to consider the second case; by applying again proposition \ref{infsum}, we obtain as in the proof of the claim an equality such as:
\[\sum_{j=1}^{s''}(\beta'_j,\varepsilon_{\beta'_j})+(0,(s_0-s')-s'')=(0,c_0-c'),\]
with $t_0-t'=s_0-s'\leq s''+(s_0-s'-s'')$. We conclude by the same reasoning as in the proof of the claim that it is incompatible with the minimality of $c_0$. $\Box$

Let $z(\Psi)$ be the quotient $\frac {c_0}{t_0}=\frac {d_0}{s_0}$; as a consequence of the claim, it doesn't depend on the choice of the $\alpha_i$ and $\beta_j$.

Now we'll return to the general case. We'll show the following result, which will imply the proposition: for every $(\alpha_1,\dots,\alpha_t)$ and $(\beta_1,\dots,\beta_s)$, we have $\frac ct\leq z(\Psi)\leq\frac ds$, and both inequalities are equalities only if $(\alpha_i,\beta_j)=0$ for every $i,j$.

Assume first $c$ and $d$ are multiples of $c_0$, say $c=c_1c_0$ and $d=d_1c_0$, and $\frac ct\geq\frac ds$; we'll show we then have $\frac ct=z(\Psi)=\frac ds$, and for every $i,j$, $(\alpha_i,\beta_j)=0$. We'll proceed by induction on $c_1+d_1$.

First remark that if $\Psi$ is $\Delta'$-complete, $-(\Phi-\Psi)$ is $\Delta'$-complete too; moreover, by replacing $\Psi$ by $-(\Phi-\Psi)$, we replace $z(\Psi)$ by $1-z(\Psi)$, $c_0$ by $t_0-c_0$ and the $\alpha_i$ (resp. the $\beta_j$) by the $-\beta_j$ (resp. the $-\alpha_i$), hence if the result is true for some $c_1,d_1$ and for $\Psi$ and $-(\Phi-\Psi)$, it is also true with $c_1$ and $d_1$ switched. Hence we may assume $\frac ct\geq z(\Psi)$, if not, then $\frac ds\leq z(\Psi)$ and we fall into the symmetrical case.

The case $c_1=d_1=1$ is simply the claim. If $c_1=1$ and $d_1>1$, then the same claim asserts $\frac ct=z(\Psi)$, hence $\frac ds\leq z(\Psi)$ and we fall again into the symmetrical case; assume then $c_1>1$.

Remark: when we'll apply the induction hypothesis, it will always be to some $\alpha'_1,\dots,\alpha'_{t'}$ and to $\beta_1,\dots,\beta_s$, hence the $\beta_j$ won't be mentioned. Moreover, since no other elements of $\Phi-\Psi$ than $\beta_1,\dots,\beta_s$ will occur in the proof, we can reduce ourselves to the case where $\Psi$ is the largest $\Delta'$-complete subset of $\Phi$ not containing them, i.e. that it contains every element of $\Phi$ which is not greater than any $\beta_j$.

Let $t'\in\{1,\dots,t\}$ be such that $\sum_{i=1}^{t'}(\alpha_i,\varepsilon_{\alpha_i})\geq(0,c_0)$ and $\sum_{i=1}^{t'-1}(\alpha_i,\varepsilon_{\alpha_i})<(0,c_0)$. By applying proposition \ref{infsum} and the induction hypothesis to $\sum_{i=1}^{t'}(\alpha_i,\varepsilon_{\alpha_i})$ (resp. $\sum_{i=t'}^t(\alpha_i,\varepsilon_{\alpha_i})$), we obtain $t'\geq \frac{c_0}{z(\Psi)})$ and $t-t'+1\geq \frac{c-c_0}{z(\Psi)}$). Since $t\leq\frac c{z(\Psi)}$, at least one of the above inequalities is an equality, and both are equalities if $t<\frac c{z(\Psi)}$.

Assume for exemple $t'=\frac{c_0}{z(\Psi)}$, the case $t-t'+1=\frac{c-c_0}{z(\Psi)}$ being symmetrical; we will then show that $(\alpha_i,\beta_j)=0$ for every $i,j$.

Set $\gamma=\sum_{i=1}^{t'}\alpha_i$. If $\gamma=0$, the assertion follows immediately from the induction hypothesis applied first to $(\alpha_1,\dots,\alpha_{t'})$ and then to $(\alpha_{t'-1},\dots,\alpha_t)$;
assume now $\gamma\neq 0$. The element $\gamma$ cannot be equal to any $\alpha_i$, since we would then have $\sum_{i'\neq i}(\alpha_{i'},\varepsilon_{\alpha_{i'}})=(0,c_0)$, which is impossible by the induction hypothesis because $\frac{c_0}{t'-1}>z(\Psi)$; since $\Phi$ is simply-laced, we then have $(\gamma,\gamma)=2$ and $(\alpha_i,\gamma)\leq 1$ for every $i$, hence there exist at least two different $i$ such that $(\alpha_i,\gamma)>0$; let's call them $i$ and $i'$. The character $\alpha_i-\gamma$ is then an element of $\Phi$, and $\alpha_i\geq\gamma$ since the equality $\sum_{i''\neq i}(\alpha_{i''},\varepsilon_{\alpha_{i''}})=(\gamma-\alpha_i,c_0+\varepsilon_{\gamma-\alpha_i})$ is impossible (apply proposition \ref{infsum} and the induction hypothesis to check this); we then have:
\[\sum_{i''\neq i}(\alpha_{i''},\varepsilon_{\alpha_{i''}})+(\alpha_i-\gamma,\varepsilon_{\alpha_i-\gamma})=(0,c_0).\]
Since $t'=\frac{c_0}{z(\Psi)}$ and $\alpha_i-\gamma\in\Psi$, we deduce from the induction hypothesis that $(\alpha_{i''},\beta_j)=0$ for every $i''\neq i$ and every $j$; by replacing $i$ by $i'$, we see that it is also true for $i''=i$.

We deduce from this that $(\gamma,\beta_j)=0$ for every $j$; with the help of the proposition \ref{infsum}, we can even replace $\gamma$ by any smaller element, which way we see that for every element $\delta$ of $\Delta'$ occurring in the decomposition of $\gamma$, $(\delta,\beta_j)=0$ for every $j$.

Moreover, we have $\gamma'=\gamma+\alpha_{t'+1}\in\Phi$; we will show that the above remarks imply $\gamma'\in\Psi$. With the proposition \ref{decprg} and an easy induction, we see that we can assume $\gamma\in\Delta'$. Suppose $\gamma'\not\in\Psi$; by the assumption made on $\Psi$, we then have $\gamma'\geq\beta_j$ for some $j$. Write:
\[\gamma'=\sum_{\delta\in\Delta'}c_\delta\delta;\]
\[\beta_j=\sum_{\delta\in\Delta'}b_\delta\delta.\]
Then $c_\delta\geq b_\delta$ for every $\delta\in\Delta'$, and $c_\gamma=b_\gamma$ since $\alpha_{t'+1}$ isn't greater than $\beta_j$. Hence:
\[(\gamma,\beta_j)=2b_\gamma-\sum_{\delta,(\delta,\gamma)=-1}b_\delta\]
\[\geq 2c_\gamma-\sum_{\delta,(\delta,\gamma)=-1}c_\delta=(\gamma,\gamma')>0,\]
which leads to a contradicton.

Consider now the equality:
\[(\gamma',\varepsilon_{\gamma'})+\sum_{i=t'+2}^t(\alpha_i,\varepsilon_{\alpha_i})=(0,c-c_0).\]
By the induction hypothesis, we have $t-t'=\frac{c-c_0}{z(\Psi)}$, hence $t=\frac c{z(\Psi)}$, and for every $j$, $(\gamma',\beta_j)=0$ and $(\alpha_i,\beta_j)=0$ for every $i\geq t'+2$; since we already know this is also true for $\alpha_i$, $i\leq t$, and $\gamma$, we obtain the desired assertion.

Assume finally $c$ is not a multiple of $c_0$; by considering the equality:
\[c_0(\sum_{i=1}^t(\alpha_i,\varepsilon_{\alpha_i}))=(0,c_0c),\]
we obtain $\frac{c_0c}{c_0t}\leq z(\Psi)$, hence $t\geq cz(\Psi)$. Since $cz(\Psi)$ is not an integer, we see the inequality is always strict. The case where $d$ is not a multiple of $d_0$ is treated similarly, and concludes the proof of the proposition. $\Box$

For every $\alpha_1,\dots,\alpha_t\in\Psi$ (resp. $\beta_1,\dots,\beta_s\in\Phi-\Psi$) whose sum is zero, if $t>0$ (resp. $s>0$), set:
\[z(\alpha_1,\dots,\alpha_t)=\frac 1t\sum_{i=1}^t\varepsilon_{\alpha_i},\]
and define $z(\beta_1,\dots,\beta_s)$ similarly.
Write:
\[z(\Psi)=\operatorname{Sup}_{(\alpha_1,\dots,\alpha_t)}z(\alpha_1,\dots,\alpha_t),\]
\[z'(\Psi)=\operatorname{Inf}_{(\beta_1,\dots,\beta_s)}z(\beta_1,\dots,\beta_s),\]
the upper (resp. lower) bound being taken over all the families $(\alpha_1,\dots,\alpha_t)$, $t>0$ (resp. $(\beta_1,\dots,\beta_s)$, $s>0$) of elements of $\Psi$ (resp. $\Phi-\Psi$) whose sum is zero; set $z(\Psi)=0$ (resp. $z'(\Psi)=1$) if there is no such family in $\Psi$ (resp. $\Phi-\Psi$). This definition of $z(\Psi)$ is clearly consistent with the one used in the proof of the proposition, and we deduce from this same proposition that $z(\Psi)\leq z'(\Psi)$.

Moreover, we have the following result:

\begin{lemme}\label{atb}
There exist $\alpha_1,\dots,\alpha_t\in\Psi$ (resp. $\beta_1,\dots,\beta_s\in\Phi-\Psi$) such that $z(\alpha_1,\dots,\alpha_t)=z(\Psi)$ (resp. $z(\beta_1,\dots,\beta_s)=z'(\Psi)$.
\end{lemme}

We will show the result for $z(\Psi)$, the proof for $z'(\Psi)$ being similar. To prove the desired assertion, we only have to show that the upper bound may be taken on a finite number of families. Let's show the following lemmas:

\begin{lemme}
Assume there exists $t'<t$ such that $\alpha_1+\dots+\alpha_{t'}=0$. Then we have:
\[z(\alpha_1,\dots,\alpha_t)\leq \operatorname{Sup}(z(\alpha_1,\dots,\alpha_{t'}),z(\alpha_{t'+1},\dots,\alpha_t)).\]
\end{lemme}

We have $z(\alpha_1,\dots,\alpha_t)=\frac{t'z(\alpha_1,\dots,\alpha_{t'})+(t-t')z(\alpha_{t'+1},\dots,\alpha_t)}t$, hence the result. $\Box$

\begin{lemme}
Assume $t>\operatorname{card}(\Phi)+1$. Then after eventually rearranging the $\alpha_i$, there exists $t'<t$ such that $\alpha_1+\dots+\alpha_{t'}=0$.
\end{lemme}

We can assume $\gamma_{t'}=\sum_{i=1}^{t'}\alpha_i$ is an element of $\overline{\Phi}$ for every $t'<t$. If one of them is zero there is nothing to prove; if all of them are nonzero, since $t>\operatorname{card}(\Phi)+1$, there exist $t'<t''$ such that $\gamma_{t'}=\gamma_{t''}$; we then have $\sum_{i=t'-1}^{t''}\gamma_i=0$, hence the result. $\Box$

According to these two lemmas, we only have to take the upper bound on the set of families $(\alpha_1,\dots,\alpha_t)$ such that $t\leq \operatorname{card}(\Phi)+1$; since this set is obviously finite, the lemma \ref{atb} is proved. $\Box$

Now we'll be concerned about the element $x$ of $\mathcal{A}$ mentioned in the theorem. Let $A_I$ be the facet of $\mathcal{B}$ associated to $I$, and $\overline{A_I}$ its closure; we have:

\begin{propo}\label{ccz}
Let $z$ be any element of $[z(\Psi),z'(\Psi)]$. There exists an element $x$ of $\overline{A_I}$ such that:
\begin{itemize}
\item for every $\alpha\in\Psi$, $\alpha(x)\leq z-\varepsilon_\alpha$;
\item for every $\alpha\in\Phi-\Psi$, $\alpha(x)\geq z-\varepsilon_\alpha$.
\end{itemize}
Moreover, if $z(\Psi)<z'(\Psi)$, the set of such elements contains an open subset of $\mathcal{A}$.
\end{propo}

Let $E_{\Psi,z}$ be the subset of the elements $x$ of $\mathcal{A}$ such that:
\begin{itemize}
\item for every $\alpha\in\Delta'$, $\alpha(x)\geq -\varepsilon_\alpha$;
\item for every maximal element $\alpha$ of $\Psi$, $\alpha(x)\leq z-\varepsilon_\alpha$;
\item for every minimal element $\alpha$ of $\Phi-\Psi$, $\alpha(x)\geq z-\varepsilon_\alpha$.
\end{itemize}
(The first condition simply amounts to say that $x\in\overline{A_I}$.)
Let $x$ be any element of $E_{\Psi,z}$; since for every $\alpha\in\Phi$ and every $\beta\in\Delta'$ such that $\alpha+\beta\in\Phi$, we have $\varepsilon_{\alpha+\beta}=\varepsilon_\alpha+\varepsilon_\beta$, we obtain by an easy induction that for every $\alpha\in\Phi$, $\alpha(x)\geq-\varepsilon_\alpha$; by a similar argument, for every $\alpha\in\Psi$, $\alpha(x)\leq z-\varepsilon_\alpha$, and for every $\alpha\in\Phi-\Psi$, $\alpha(x)\geq z-\varepsilon_\alpha$; hence $x$ satisfies the condition of the lemma. We then only have to show that $E_{\Psi,z}$ is nonempty. 

Let $\Phi_m$ be the subset of $\Phi$ containing $-\Delta'$, the maximal elements of $\Psi$ and the opposites of the mimimal elements of $\Phi-\Psi$; set, for every $\alpha\in\Phi_m$:
\begin{itemize}
\item if $\alpha\in\Psi$, $f_z(\alpha)=z-\varepsilon_\alpha$;
\item if $\alpha$ is the opposite of some minimal element of $\Phi-\Psi$, $f_z(\alpha)=-z+\varepsilon_{-\alpha}=1-z-\varepsilon_\alpha$;
\item if $\alpha\in-\Delta'$, $f_z(\alpha)=\varepsilon_{-\alpha}=1-\varepsilon_\alpha$;
\item if $\alpha$ satisfies more than one of the above conditions, $f_z$ takes the lowest possible value.
\end{itemize}

The set $E_{\Psi,z}$ is nonempty if and only if for every $\alpha_1,\dots,\alpha_t\in\Phi_m$ whose sum is zero, $\sum_{i=1}^tf_z(\alpha_i)\geq 0$. The first assertion of the proposition follows then from the following result:

\begin{lemme}\label{ccz1}
Let $\alpha_1,\dots,\alpha_t$ (resp. $\beta_1,\dots,\beta_s$) be elements of $\Psi$ (resp. $\Phi-\Psi$) such that $\sum_{i=1}^t\alpha_i=\sum_{j=1}^s\beta_j$. Set:
\[c_0=\sum_{i=1}^t\varepsilon_{\alpha_i}+\sum_{j=1}^s\varepsilon_{-\beta_j}.\]
Then $tz+s(1-z)\geq c_0$, and the inequality is strict if $z(\Psi)<z'(\Psi)$ and $s,t>0$.
\end{lemme}

Write $\gamma=\sum_{i=1}^t\alpha_i=\sum_{j=1}^t\beta_j$. Note that $\gamma$ is not necessarily an element of $\overline{\Phi}$.

We'll show the lemma by induction on $s$. Assume first $\gamma=0$, and let $c'_0$ (resp. $c''_0$) be the sum of the $\varepsilon_{\alpha_i}$ (resp. the $\varepsilon_{\beta_j}$). It follows immediately from the definition of $z(\Psi)$ and $z'(\Psi)$ that we have $tz\geq c'_0$ and $sz\leq c''_0$, and that at least one of these inequalities is strict if $z(\Psi)<z'(\Psi)$ and $s,t>0$; hence:
\[tz+s(1-z)\geq c'_0+s-c''_0=c_0,\]
and the inequality is strict if $z(\Psi)<z'(\Psi)$ and $s,t>0$.

Moreover, since for each $j$, $\varepsilon_{-\beta_j}=1-\varepsilon_{\beta_j}$, we obtain:
\[c'_0+s-c''_0=\sum_{i=1}^t\varepsilon_{\alpha_i}+\sum_{j=1}^s\varepsilon_{-\beta_j}=c_0,\]
which proves $tz+s(1-z)\geq c_0$ ($>c_0$ if $z(\Psi)<z'(\Psi)$ and $s,t>0$). Assume now $\gamma\neq 0$; we now have:
\[(\sum_{i=1}^t\alpha_i,\sum_{j=1}^s\beta_j)>0;\]
there exist then $i_0\in\{1,\dots,t\}$ and $j_0\in\{1,\dots,s\}$ such that $(\alpha_{i_0},\beta_{j_0})>0$; we will assume $i_0=t$ and $j_0=s$. We have $\alpha_t-\beta_s\in\Phi$, and since $\alpha_t\in\Psi$ and $\beta_s\not\in\Psi$, $\beta_s=\alpha_t+(\beta_s-\alpha_t)>\alpha_t,\beta_s-\alpha_t$, hence $\alpha_t-\beta_s>-\beta_s,\alpha$; we then have $h(\alpha_t)+h(-\beta_s)=h(\alpha_t-\beta_s)<h(\Phi)$. Consider now the equality:
\[\sum_{i=1}^{t-1}(\alpha_i,\varepsilon_{\alpha_i}+\sum_{j=1}^{s-1}(-\beta_j,\varepsilon_{-\beta_j})=(-\alpha_t+\beta_s,c_0-1).\]
According to the proposition \ref{infsum}, there exist $\alpha'_1,\dots,\alpha'_t\in\Psi\cup\{0\}$ (resp. $\beta'_1,\dots,\beta'_s\in(\Phi-\Psi)\cup\{0\}$) such that for every $i$ (resp. $j$), $\alpha'_i\leq\alpha'_j$ (resp. $\beta'_i\geq\beta'_j$) and:
\[\sum_{i=1}^{t-1}(\alpha'_i,\varepsilon_{\alpha'_i})+\sum_{j=1}^{s-1}(-\beta_j,\varepsilon_{-\beta_j})=(0,c_0-1),\]
and we deduce from the induction hypothesis:
\[(s-1)z+(t-1)(1-z)\geq c_0-1,\]
witn a strict inequality if $z(\Psi)<z'(\Psi)$ and $s,t>0$; hence the result. $\Box$

Assume now $z(\Psi)<z'(\Psi)$; we'll show that $E_{\Psi_z}$ contains an open subset of $\mathcal{A}$.
Assume this is not the case; since $E_{\Psi_z}$ is obviously convex, it is then contained in the hyperplane defined by some equation $\gamma(y)=\lambda$, with $\gamma\in X^*(T)$ and $\lambda\in{\mth{R}}$, hence:
\begin{itemize}
\item there exist $\alpha'_1,\dots,\alpha'_{t'}\in\Phi_m$ whose sum is $\gamma$ and such that $\sum_{i=1}^{t'}f_z(\alpha'_i)=\lambda$;
\item there exist $\beta'_1,\dots,\beta'_{s'}\in\Phi_m$ whose sum is $-\gamma$ and such that $\sum_{j=1}^{s'}f_z(\beta'_j)=-\lambda$.
\end{itemize}

We then have $\sum_{i=1}^{t'}\alpha'_i+\sum_{j=1}^{s'}\beta'_j=0$ and $\sum_{i=1}^{t'}f_z(\alpha'_i)+\sum_{j=1}^{s'}f_z(\beta'_j)=0$; this is possible only if all the inequalities occuring in the proof of the previous lemma, when applied to those elements, are equalities, and in particular if $tz=c'_0$ and $sz=c''_0$, which implies $z(\Psi)=z=z'(\Psi)$. Hence the result. $\Box$

\section{Proof of the theorem \ref{cln}}

In this section, we'll prove the main theorem. Let's begin by the following preliminary result:

\begin{lemme}
Assume $H=T'U_{x,r}$ for some $T',x,r$. Then $H\subset G_{x,r}$.
\end{lemme}

Since $T'$ fixes $\mathcal{A}$ pointwise, $H$ fixes $B(x,r)\cap\mathcal{A}$ pointwise. Moreover, let $y$ be any element of $B(x,r)$, and let $\mathcal{A}'$ be an apartment of $\mathcal{B}$ containing $A_I$ and $y$ (such an apartment exists by \cite[I.2.3.1]{bt}). According to \cite[I.2.5.8]{bt}, there exists $g\in G^0$ such that $g(y)\in\mathcal{A}$ and $g$ fixes $A_I$ pointwise, hence $g\in I$; $y$ is then fixed by $gHg^{-1}=H$. Hence $H$ fixes $B(x,r)$ pointwise, which amounts to say it is contained in $G_{x,r}$, as required. $\Box$

Now let's go on into the proof of the theorem. First we'll suppose $\Phi$ is of type $1$; let $\Psi$ be the subset of the elements $\alpha\in\Phi$ such that $f_\varepsilon(\alpha)=v'+1$. We have already seen $\Psi$ is $\Delta'$-complete; let's show the following result:

\begin{propo}\label{cnc}
Let $z$ be any element of $[z(\Psi),z'(\Psi)]$; set $r=v'+z$. Let $x$ be any element of $E_{\Psi,z}$; we have:
\[U_{x,r}^+\subset H\subset G_{x,r}.\]
\end{propo}

Let $f_{x,r}$ be the concave function associated to $G_{x,r}$; we have, for every $\alpha\in\Phi$:
\[f_{x,r}(\alpha)=\operatorname{ceil}(-\alpha(x)+r).\]
With the help of the previous lemma, in order to show that $H\subset G_{x,r}$, we only have to prove that $f\geq f_{x,r}$. Let $\alpha$ be any element of $\Phi$; if $\alpha\in\Psi$, we have, using the definition of $E_{\Psi,z}$ and the fact that $f_H(\alpha)$ is an integer:
\[f_H(\alpha)=v'+1+\varepsilon_\alpha=r+1-(-\varepsilon_\alpha)-z\]
\[\geq \operatorname{ceil}(-\alpha(x)+r-z)\geq f_{x,r}(\alpha),\]
and if $\alpha\not\in\Psi$:
\[f_H(\alpha)=v'+\varepsilon_\alpha=r-(z-\varepsilon_\alpha)\geq \operatorname{ceil}(-\alpha(x)+r),\]
hence the result.

Let now $f_{x,r}^+$ be the concave function associated to $G_{x,r}^+$; we have for every $\alpha\in\Phi$:
\[f_{x,r}^+(\alpha)=\operatorname{Inf}_{r'>r}f_{\alpha,x,r'}(\alpha)=\operatorname{ceil}(-\alpha(x)+r')\]
\[=\operatorname{floor}(-\alpha(x)+r+1);\]
we'll show that $f\leq f_{x,r}^+$, which will imply $U_{x,r}^+\subset H$. For every $\alpha\in\Phi$, we have, when $\alpha\in\Psi$:
\[f_H(\alpha)=r+1-(z-\varepsilon_\alpha)\leq \operatorname{floor}(r+1-\alpha(x)),\]
and when $\alpha\not\in\Psi$:
\[f_H(\alpha)=r+1-(1-\varepsilon_\alpha)-z\]
\[\leq \operatorname{floor}(r+1-\alpha(x)-z)\leq f_{x,r}^+(\alpha),\]
hence again the result. $\Box$

First consider the case $z(\Psi)<z'(\Psi)$; we will show that the groups $H$ such that $U_{x,r}\subset H\subset G_{x,r}$ for some $x,r$ are the ones which satisfy this condition. Since we can choose $z=r-v'$ anywhere in $[z(\Psi),z'(\Psi)]$, there exist $r<r'$ and $x,x'$ such that the previous proposition holds for both $x,r$ and $x',r'$. In particular we obtain:
\[U_{x,r''}\subset H\subset G_{x',r'}\]
for any $r''>r$. Unfortunately this isn't enough since the inclusion $U_{x,r''}\subset U_{x',r'}$ may very well be strict even if $r''<r'$; this problem will be solved by checking that we may assume $x=x'$, which is done by the following lemma:

\begin{lemme}
There exist $z<z'\in[z(\Psi),z'(\Psi)]$ such that $E_{\Psi,z}\cap E_{\Psi,z'}\neq\emptyset$.
\end{lemme}

Since $z(\Psi)<z'(\Psi)$, for every $z$, $E_{\Psi_z}$ contains an open subset of $A_I$, hence its volume is nonzero. Assume all $E_{\Psi_z}$ are disjoint; then $A_I$ contains an uncountable union of disjoint subsets with positive volume, which is impossible. $\Box$

\begin{coro}\label{cas}
Set $r'=z'+v'$; we have $U_{x,r'}\subset H\subset G_{x,r'}$.
\end{coro}

According to the proposition \ref{cnc} and the previous lemma, we now have $U_{x,r'}\subset U_{x,r}^+\subset H\subset G_{x,r'}$. $\Box$

Hence the theorem is proved for $f_H$ of type $1$ and such that $z(\Psi)<z'(\Psi)$. In this case, we may show the following converse:

\begin{propo}
Assume $U_{x,r}\subset H\subset G_{x,r}$ for some $x,r$, and set $r=z+v'$, with $z\in[0,1]$ and $v'$ being an integer. Then $f_H$ is of type $1$ and $z(\Psi)<z'(\Psi)$.
\end{propo}

For every $\alpha\in\Phi$, we have:
\[f_\varepsilon(\alpha)=\operatorname{ceil}(-\alpha(x)+r)-\varepsilon_\alpha\in\{\operatorname{ceil}(r-1),\operatorname{ceil}(r)\},\]
hence $f_H$ is of type $1$. Moreover, since the valuation is discrete, we have $H=G_{x,r}$ for every $r$ in some interval $]r',r'']$, with $r'<r''$; we obviously can assume $\operatorname{floor}(r')=\operatorname{floor}(r'')=v'$. Since $f_H(x)=\operatorname{ceil}(-\alpha(x)+r)$ for every $\alpha\in\Phi$ and every $r\in]r',r'']$, setting $z'=r'-v'$ and $z''=r''-v'$, we obtain that every element $\alpha$ of $\Psi$ (resp. $\Phi-\Psi$) must satisfy $\alpha(x)<z'$ (resp. $\alpha(x)\geq z''$), hence $z(\Psi)\leq z'<z''\leq z'(\Psi)$. $\Box$

Suppose now $f_H$ is still of type $1$, but such that $z(\Psi)=z'(\Psi)$; according to the previous proposition, $H$ cannot be equal to any group of the form $T'U_{x,r}$. Such a case actually occurs: for example, take $G=SL_4$, and:
\[H=\left\{\left(\begin{array}{cccc}1+{\mathfrak{p}}&{\mathfrak{p}}&{\mathfrak{p}}&\mathcal{O}\\{\mathfrak{p}}&1+{\mathfrak{p}}&{\mathfrak{p}}&\mathcal{O}\\{\mathfrak{p}}^2&{\mathfrak{p}}&1+{\mathfrak{p}}&{\mathfrak{p}}\\{\mathfrak{p}}^2&{\mathfrak{p}}&{\mathfrak{p}}&1+{\mathfrak{p}}\end{array}\right)\right\}.\]
It is easy to check that this group is a normal subgroup of $I$ and that the corresponding concave function $f_H$ is of type $1$. Moreover, we have:
\[\Psi=\{\alpha_1,\alpha_2,\alpha_1+\alpha_2,\alpha_3,-\alpha_M,\alpha_3-\alpha_M.\}\]
Since $(\alpha_1+\alpha_2)+(\alpha_3-\alpha_M)=0$, we have $z(\Psi)\geq \frac 12$, and since $(\alpha_2+\alpha_3)+(\alpha_1-\alpha_M)=0$, we have $z'(\Psi)\leq\frac 12$. Hence $z(\Psi)=z'(\Psi)=\frac 12$.

Consider the families $(\beta_1,\dots,\beta_s)$ of elements of $\Phi-\Psi$ such that $\sum_{j=1}^s(\beta_j,\varepsilon_{\beta_j})=(0,d)$ with $\frac ds=z'(\Psi)$; let's call them $z'(\Psi)$-families. Let $\Phi_0$ be the root subsystem of elements of $\Phi$ which are linear combinations of elements of such families; we have:

\begin{propo}\label{corth}
The rank of the root subsystem $\Phi_0$ of $\Phi$ is strictly smaller than the rank of $\Phi$.
\end{propo}

Since $z(\Psi)=z'(\Psi)$, there exist $\alpha_1,\dots,\alpha_t\in\Psi$ and an integer $c$ such that $\sum_{i=1}^t(\alpha_i,\varepsilon_{\alpha_i})=(0,c)$ and $\frac ct=z(\Psi)$. According to the proposition \ref{smph}, every element of $\Phi_0$ is then strongly orthogonal to all the $\alpha_i$, and this implies the result. $\Box$

Moreover, set $\Psi_0=\Psi\cap\Phi_0$; we obviously have $z'(\Psi_0)=z'(\Psi)$, and since the $\alpha_i$ of the previous proposition cannot belong to $\Phi_0$, $z(\Psi_0)$ is strictly smaller than $z(\Psi)$. Let $M$ be the Levi subgroup of $G$ generated by $T$ and the $U_\alpha$, $\alpha\in\Phi_0$; $M$ is proper, and since $H\cap M$ is obviously normal in the Iwahori subgroup $I\cap M$ of $M$, according to the corollary \ref{cas}, there exists an element $x$ of the apartment $\mathcal{A}_M$ of the Bruhat-Tits building of $M$ associated to $S$ and an element $r$ of $[v',v'+1]\subset{\mth{R}}^+$ such that $U_{M,x,r}\subset H\cap M\subset M_{x,r}$; moreover, using the canonical projection $\mathcal{A}\rightarrow\mathcal{A}_M$, with a slight abuse of notation, setting $z=r-v'$, we will also call $E_{\Psi_0,z}$ the convex subset of $\mathcal{A}$ whose image in $\mathcal{M}$ is $E_{\Psi_0,z}$; we may then consider $x$ as an element of $\mathcal{A}$.

Now we'll determine $r'$. we have the following result:

\begin{lemme}\label{psi1}
Let $\Psi'$ be the subset of $\Phi$ which is the union of $\Psi$ and all $z'(\Psi)$-families. Then $\Psi'$ is $\Delta'$-complete. and $z(\Psi)=z(\Psi')$.
\end{lemme}

Let's show the first assertion: we will in fact show the equivalent assertion that $-(\Phi-\Psi')$ is $\Delta'$-complete. Let $\alpha$ be any element of $\Phi-\Psi'$. Assume there exists $\beta_1,\dots,\beta_t\in\Phi-\Psi$ such that $\sum_{i=1}^t(\beta_i,\varepsilon_{\beta_i})=(0,d)$, with $d=tz'(\Psi)$, and $\beta_t\geq\alpha$; we then have:
\[\sum_{i=1}^{t-1}(\beta_i,\varepsilon_{\beta_i})\leq (-\alpha,d-1+\varepsilon_{-\alpha}).\]
According to the proposition \ref{infsum}, there exist then $\beta'_1,\dots,\beta'_{t-1}\in\Phi-\Psi$ (one may check as in the proof of the proposition \ref{smph} that assuming that at least one of the $(\beta_i,\varepsilon_{\beta_i})$ is reduced to $(0,1)$ leads to a contradiction) such that:
\[\sum_{i=1}^{t-1}(\beta'_i,\varepsilon_{\beta'_i})=(-\alpha,d-1+\varepsilon_{-\alpha}),\]
hence:
\[\sum_{i=1}^{t-1}(\beta'_i,\varepsilon_{\beta'_i})+(\alpha,\varepsilon_\alpha)=(0,d).\]
Hence the family $(\beta'_1,\dots,\beta'_{t-1},\alpha)$ is a $z'(\Psi)$-family. Since $\alpha\not\in\Psi'$, this is impossible; $\beta_t$ cannot then be greater than $\alpha$. Since we already know that no element of $\Psi$ can be greater than $\alpha$, we have just shown that every element of $\Phi$ which is greater than $\alpha$ belongs to $\Phi-\Psi'$; hence $-(\Phi-\Psi')$ is $\Delta'$-complete.

Now we'll prove that $z(\Psi)=z(\Psi')$. Let $\alpha_1,\dots,\alpha_t$ be elements of $\Psi'$ such that $\sum_{i=1}^t(\alpha_i,\varepsilon_{\alpha_i})=(0,c)$, with $\frac ct\geq z(\Psi)$. Assume $\alpha_i\in\Psi$ for every $i\leq t'$, and $\alpha_i$ belongs to some $z(\Psi)$-family if $i>t'$; in the second case, let $\beta_{i,1},\dots,\beta_{i,s_i-1}$ be the other members of that family. Setting $d_i=\frac{s_i}{z(\Phi)}$ for every $i>t'$, we obtain:
\[\sum_{i=1}^{t'}(\alpha_i,\varepsilon_{\alpha_i})+\sum_{i=t'+1}^t(0,d_i)=\sum_{i=t'+1}^t\sum_{j=1}^{s_i-1}(\beta_{i,j},\varepsilon_{\beta_{i,j}})+(0,c).\]
By the same process as in the proof of the proposition \ref{smph}, we can obtain an equality:
\[0=\sum_{i=1}^{t''}(\alpha'_i,\varepsilon_{\alpha'_i})+\sum_{i=t'+1}^t(0,d_i)-\sum_{j=1}^{s'}(\beta'_j,\varepsilon_{\beta'_,j})\]
\[+(0,-c-\sum_{i=t'+1}^t(s_i-1)+s'),\]
where $t''\leq t'$, $s'\leq\sum_{i=t'+1}^t(s_i-1)$ and the $\alpha'_i$ (resp. the $\beta'_j$) are elements of $\Psi$ (resp. $\Phi-\Psi)$ such that $(\alpha'_i,\beta'_j)=0$ for every $i,j$. The sum of the $\alpha'_i$ (resp. $\beta'_j$) must then be zero; writing $(0,c')=\sum_{i=1}^{t''}(\alpha'_i,\varepsilon_{\alpha'_i})$ (resp. $(0,d')=\sum_{j=1}^{s'}(\beta'_j,\varepsilon_{\beta'_,j})$), we obtain:
\[c'+\sum_{i=t'+1}^td_i-d'-c-\sum_{i=t'+1}^t(s_i-1)+s'=0,\]
hence, since by definition $c'\leq t''z(\Psi)$ and $d'\geq s'z'(\Psi)=s'z(\Psi)$, and by assumption $c\geq tz(\Psi)$:
\[(t''+\sum_{i=t'+1}^ts_i-s'-t)z(\Psi)-\sum_{i=t'+1}^t(s_i-1)+s\geq 0,\]
which can be rewritten as:
\[(t''-t')z(\Psi)+(\sum_{i=t'+1}^t(s_i-1)-s')(z(\Psi)-1)\geq 0.\]
Since $0<z(\Psi)<1$, $t''\leq t'$ and $\sum_{i=t'+1}^t(s_i-1)\geq s'$, this is possible only if all the inequalities we have combined to get the above one are equalities, and in particular if $v=tz(\Psi)$. Hence $z(\Psi')\leq z(\Psi)$; since the other inequality is obvious, the lemma is proved. $\Box$

According to the first lemma, if $f'_H$ is the concave function on $\Phi$ associated to $\Psi_1$, $U_{f'_H}$ is normal in $I$; moreover, we obviously have $z'(\Psi')>z'(\Psi)$ hence $z(\Psi')<z'(\Psi')$ by the second lemma. There exist then $x'\in\mathcal{A}$ and $r'>r\in{\mth{R}}$ such that $U_{f'_H}=U_{x',r'}$.

Moreover, since $H\cap M$ is contained in $I$, it normalizes $U_{f'_H}$. Since for every $\alpha\in\Phi$, we have $f'_H(\alpha)\geq f_H(\alpha)$ if $\alpha\in\Psi_0$ and $f'_H(\alpha)=f_H(\alpha)$ if $\alpha\not\in\Psi_0$, we obtain:
\[H=(H\cap M)U_{x',r'}=T'U_{M,x,r}U_{x',r'}\]
for some $T'\subset P_T$.

It remains to see that we can choose $x=x'$, or equivalently, that the intersection $E_{\Psi_0,z}\cap E_{\Psi',z'}$ is nonempty for some $z'>z$. Let $\Psi''$ be the smallest $\Delta'$-complete subset of $\Phi$ containing $\Psi_0$; we have $z(\Psi'')<z(\Psi)$, and $z'(\Psi'')=z'(\Psi)$ by a similar argument as in the lemma \ref{psi1}. Moreover, we have:

\begin{lemme}\label{ccz3}
Let $\alpha_1,\dots,\alpha_t$ (resp. $\beta_1,\dots,\beta_s$, $\gamma_1,\dots,\gamma_r$,$\delta_1,\dots,\delta_q$) be elements of $\Psi''$ (resp. $\Phi-\Psi''$, $\Psi'$, $\Phi-\Psi'$) such that $\sum_{i=1}^t\alpha_i+\sum_{k=1}^r\gamma_k=\sum_{j=1}^s\beta_j+\sum_{l=1}^q\delta_l$. Set:
\[c_0=\sum_{i=1}^t\varepsilon_{\alpha_i}+\sum_{j=1}^s\varepsilon_{-\beta_j}+\sum_{k=1}^r\varepsilon_{\gamma_k}+\sum_{l=1}^q\varepsilon_{-\delta_l}.\]
Then $tz+s(1-z)+rz'+q(1-z')\geq c_0$ for $z'>z(\Psi)>z$ and $z$ and $z'$ close enough to each other.
\end{lemme}

If there exist $i,j$ such that $(\alpha_i,\beta_j)>0$, we can use the induction hypothesis as in the proof of the lemma \ref{ccz1}; assume that $(\alpha_i,\beta_j)\leq 0$ for every $i,j$. For the same reason, we can also assume $(\gamma_k,\delta_l)\leq 0$ for every $k,l$, and since $z'+1-z\geq 1$, $(\beta_j,\gamma_k)\leq 0$ for every $j,k$; moreover, if $z'\geq 1-z$ (resp. $z'\leq 1-z$), we have $z+z'\geq 1$ (resp. $(1-z)+(1-z')\geq 1$, hence we can assume $(-\alpha_i,\gamma_k)\leq 0$ for every $i,k$ (resp. $(-\beta_j,\delta_l)\leq 0$ for every $j,l$). We will suppose $z+z'\geq 1$, the other case being similar.

Consider the equality:
\[\sum_{k=1}^r\gamma_k=\sum_{j=1}^s\beta_j+\sum_{l=1}^q\delta_l+\sum_{i=1}^t(-\alpha_i).\]
Since the product of any term of the left-hand side by any term of the right-hand side is smaller than $0$, all of them must be zero; in particular, we have $\sum_{k=1}^r\gamma_k=0$. Since we then obtain $s'z'\geq s''z(\Psi)\geq\sum_{k=1}^r\varepsilon_{\gamma_k}$, we may assume $r=0$.

For any $z\in[z(\Psi''),z(\Psi)]$, according to the lemma \ref{ccz1} applied to $\Psi''$ and $z$, we have:
\[tz+(s+r)(1-z)>c_0.\]
This is in particular true for $z=z(\Psi)$; we easily deduce from this that we have $tz+s(1-z)+r(1-z')>c_0$ for $z\leq z(\Psi)<z'$ and $z$ and $z'$ close enough to $z(\Psi)$, as required. $\Box$

By the same argument as in the proof of the lemma \ref{atb}, we only have to consider a finite number of families $((\alpha_i),(\beta_j),(\gamma_k),(\delta_l))$; there exist then $z,z'$ such that the above lemma is true for $z,z'$ and all such families. Hence $E(\Psi'',z)\cap E(\Psi',z')$ is nonempty; since $E(\Psi'',z)$ is obviously contained in $E(\Psi_0,z)$, we obtain the desired result.

We'll now turn on to the case when $f_H$ is of type $2$. This case actually occurs too: for exemple, take $G=SL_4$, and:
\[H=\left\{\left(\begin{array}{cccc}1+{\mathfrak{p}}&{\mathfrak{p}}&{\mathfrak{p}}&{\mathfrak{p}}\\{\mathfrak{p}}&1+{\mathfrak{p}}&{\mathfrak{p}}&{\mathfrak{p}}\\{\mathfrak{p}}^2&{\mathfrak{p}}^2&1+{\mathfrak{p}}^2&{\mathfrak{p}}^2\\{\mathfrak{p}}^2&{\mathfrak{p}}^2&{\mathfrak{p}}^2&1+{\mathfrak{p}}^2\end{array}\right)\right\}.\]
It is easy to check that this group is a normal subgroup of $I$. Moreover, we have $f_\varepsilon(-\alpha_1)=0$ and $f_\varepsilon(\alpha_3)=2$, hence $f_H$ is of type $2$.

Once again, $h$ is not of the form $T'U_{x,r}$ for any $T',x,r$. Let $\Phi_0$ be the subsystem of $\Phi$ whose elements are the linear combinations of elements $\alpha$ such that $f_\varepsilon(\alpha)=v'$; we have:

\begin{propo}\label{corth2}
The root subsystem $\Phi_0$ of $\Phi$ is of rank strictly smaller than the rank of $\Phi$. Moreover, it admits a subset of $\Delta'$ as a set of simple roots.
\end{propo}

We'll show that if $\alpha,\beta$ are elements of $\Phi$ such that $f_\varepsilon(\alpha)=v'$ and $f_\varepsilon(\beta)=v'+2$, $\alpha$ and $\beta$ are strongly orthogonal; we'll then obtain the first assertion of the proposition the same way as in the proposition \ref{corth}.

Assume $\alpha+\beta$ belongs to $\Phi$. If $\alpha+\beta$ is greater than $\alpha$ and $\beta$, we have:
\[f_\varepsilon(\beta)\leq f_\varepsilon(\alpha+\beta)+1\leq f_\varepsilon(\alpha)+1=v'+1,\]
which is impossible; the case $\alpha+\beta$ lesser than $\alpha$ and $\beta$ is impossible too for similar reasons. Assume now $\alpha-\beta$ belongs to $\Phi$; we have $\beta=(\alpha-\beta)+\alpha$, hence $f_\varepsilon(\beta)\leq f_\varepsilon(\alpha)+1$, which is again impossible. We then obtain that $\alpha$ and $\beta$ are strongly orthogonal.

Moreover, if $f_\varepsilon(\alpha)=0$, then $f_\varepsilon(\alpha')=0$ for every $\alpha'\geq\alpha$; we then obtain that every element of $\Delta'$ occurring in the decomposition of $-\alpha$ belongs to $\Phi_0$, which shows the second assertion. $\Box$

Consider the restriction of $f_H$ to $\Phi_0$. The subset $\Phi_0$ doesn't contain any $\alpha$ such that $f_\varepsilon(\alpha)=v'+2$, since we have just seen such an $\alpha$ is strongly orthogonal to $\Phi_0$; hence $f_H|_{\Phi_0}$ is of type $1$. Moreover, let $\Delta'_0$ be the extended basis of $\Phi_0$ whose elements are the minimal elements of $\Phi_0$, and set $\Psi_0=\Psi\cap\Phi_0$; $\Psi_0$ is obviously $\Delta'_0$-complete, and we have the following result:

\begin{lemme}
We have $z'(\Psi_0)=1$.
\end{lemme}

Consider the extended basis $\Delta'_0$ of $\Phi_0$; this is the union of $\Delta_0=\Delta'\cap\Phi_0$ with one additional element $\alpha$, which is the inverse of the greatest root of $\Phi_0$ w.r.t $\Delta_0$. Assume $f_H(-\alpha)=v'$; then for every $\beta\in\Phi$ such that $\beta\leq\alpha$, we must then have $f_H(-\beta)=v'$, hence $\beta\in\Phi_0$. By minimality of $\alpha$ in $\Phi_0$, we obtain that no such $\beta$ exists, hence $\alpha\in\Delta'$, which leads to a contradiction. Therefore $\Phi_0-\Psi_0$ doesn't contain the whole set $-\Delta'_0$, from which we deduce that $z'(\Psi_0)=1$. $\Box$

Since $z(\Psi_0)<1$ by definition, we may define $M,x$ and $r$ the same way as in the case $f_H$ of type $1$ and $z(\Psi)=z'(\Psi)$; we obtain $r\in[v',v'+1[$.

Moreover, let $f'_H$ be the function on $\Phi$ defined by $f'_H(\alpha)=f_H(\alpha)+1$ if $f_\varepsilon(\alpha)=0$, and $f'_H(\alpha)=f_H(\alpha)$ if $f_\varepsilon(\alpha)>0$; we have:

\begin{lemme}
The group $U_{f'_H}$ is normal in $I$.
\end{lemme}

Let $\alpha,\beta$ be elements of $\Phi$ such that $\alpha+\beta\in\Phi$. Suppose that we have:
\[f'_H(\alpha+\beta)>f_I(\alpha)+f'_H(\beta).\]
This is possible only if we have $f_H(\alpha+\beta)=f_I(\alpha)+f_H(\beta)$, $f'_H(\alpha)=f_H(\alpha)-1$ and $f'_H(\beta)=f_H(\beta)$. The first equality implies:
\[f_\varepsilon(\alpha+\beta)\geq f_\varepsilon(\beta).\]
Since $f_\varepsilon(\alpha+\beta)=0$ and $f_\varepsilon(\beta)>0$, this is impossible and the lemma is proved. $\Box$

It is obvious from its definition that $f'_H$ is of type $1$. Moreover, let $\Psi'$ be the $\Delta'$-complete subset of $\Phi$ associated to $f'_H$; the elements $\alpha$ of $\Delta'$ contained in $\Psi'$ are exactly the ones such that $f_\varepsilon(\alpha)=2$. Since, according to the proof of the proposition \ref{corth2}, the corresponding subset of $\Delta'$ is strongly orthogonal to some other nonempty subset of $\Phi$, $\Psi'$ cannot contain the whole set $\Delta'$; hence $z(\Psi')=0$. Since $z'(\Psi')>0$ by definition, we may define $x'$ and $r'$ as in the previous case too; we obtain here $r'\in[v'+1,v'+2[$.

We'll again conclude by proving that we may choose $x=x'$, or equivalently that $E_{\Psi_0,z}\cap E_{\Psi',z'}$ is nonempty for some $z,z'$. Let $x$ be an element of $\mathcal{A}$ satisfying the following conditions:
\begin{itemize}
\item $x\in E_{\Psi_0,z}$ for some $z\geq z(\Psi_0)$, and $\alpha(x)>0$ for every $\alpha\in\Delta'\cap\Psi_0$;
\item $x\in E_{\Psi',z'}$ for some $z'$ such that $z'<\alpha(x)$ for every $\alpha\in\Delta'\cap\Psi_0$.
\end{itemize}
The first condition is possible because $E_{\Psi_0,z}$ contains an open subset of $\mathcal{A}$, and the fact that $\Psi_0$ and $\Psi'$ are strongly orthogonal to each other allows us to add the second one; we obtain that way $U_{x,r'}U_{M,x,r}\subset H\subset G_{x,r'}M_{x,r}$, as required.

We'll finally assume $H$ is a normal subgroup of $P$, where $P$ is any parahoric subgroup of $G$ containing $I$. Since $I$ normalizes $H$, $H$ is either of the form $T'U_{x,r}$ for some $T',x,r$ or of the form $T'U_{M,x,r}U_{x,r'}$ for some $T',x,r,M,r'$; it only remains to show that we can choose $x$ in $\overline{A_P}$.

Let $P^+$ be the pro-nilpotent radical of $P$; let $P_S$ be the unique parahoric subgroup of $S(F)$ and let $P_S^+$ be its pro-nilpotent radical. For every $w$ in the Weyl group $W_P$ of $P/P^+$ relatively to $P_S/P_S^+$, $w$ normalizes $H$; hence $H$ is equal to $w(T')U_{w(x),r}$ (resp. $w(T')U_{M,w(x),r}U_{w(x),r'}$). The subset $B_H$ of $\overline{A_I}$ containing all elements $x$ satisfying the condition of the theorem is then stable by $W_P$; moreover, for any $x,x'\in B_H$ and every $t\in[0,1]$, the element $(1-t)x+tx'$ is well-defined in every apartment of $\mathcal{B}$ containing both $x$ and $x'$, and we have:
\[B(x,r)\cap B(x',r)\subset B((1-t)x+tx',r)\subset E(B(x,r),B(x',r)),\]
where $E(B(x,r),B(x',r))$ is the closure of $B(x,r)\cup B(x',r)$ in $\mathcal{B}$; hence $B((1-t)x+tx',r)$ is fixed by $H$ too. We deduce from this that $B_H$ is a convex subset of $\overline{A_I}$; since $W_P$ is finite, for any $x\in B_H$, the barycenter of the $w(x)$, $w\in W_P$, is both contained in $B_H$ and fixed by $W_P$, hence contained in $\overline{A_P}$. This completes the proof of the theorem. $\Box$

\section{The general case}

In this last section, we don't suppose $G$ to be quasi-simple anymore. Let $\Phi_1,\dots,\Phi_k$ be the connected components of $\Phi$; for every $i\in\{1,\dots,k\}$, let $G_i$ be the subgroup of $G(F)$ generated by $T(F)$ and the $U_{\alpha}(F)$, $\alpha\in\Phi_i$. Let $P$ be a parahoric subgroup containing $P_T$, and let $H$ be a normal subgroup of $P$; an immediate consequence of the proposition \ref{nsgf} is that we have:
\[H=\prod_{i=1}^k(H\cap G_i).\]
For every $i$, let $\mathcal{B}_i$ be the Bruhat-Tits building of $G_i$, and let $\mathcal{A}_i$ be the apartment of $\mathcal{B}_i$ associated to $S(F)$; it is easy to check that $\mathcal{A}$ is canonically isomorphic to $\prod_{i=1}^k\mathcal{A}_i$. Moreover, let $A_P$ be the facet of $\mathcal{A}$ attached to $P$; $A_P$ is canonically isomorphic to $\prod_{i=1}^kA_{P_i}$, where for every $i$, $A_{P_i}$ is the facet of $\mathcal{A}_i$ attached to $P_i=P\cap G_i$.

Let $x=(x_1,\dots,x_k)$ be an element of $\mathcal{A}$, and let $r_1,\dots,r_k$ be nonnegative real numbers; we'll set:
\[G_{x,(r_1,\dots,r_k)}=\prod_{i=1}^k(G_i)_{x_i,r_i},\]
where for every $i$, $(G_i)_{x_i,r_i}$ is the standard filtration subgroup of $G_i$ attached to $x_i$ and $r_i$. If $x$ belongs to $A_P$, $G_{x,(r_1,\dots,r_k)}$ is obviously normal in $P$. We'll define $U_{x_i,r_i}$ and $U_{x,(r_1,\dots,r_k)}$ in a similar fashion.

According to theorem \ref{cln}, for every $i$, there exists $x_i\in\overline{A_{P_i}}$ and $r_i\in{\mth{R}}^+$ such that $H\cap G_i$ satisfies one of the following conditions:
\begin{itemize}
\item $U_{x_i,r_i}\subset H\cap G_i\subset (G_i)_{x_i,r_i}$;
\item there exists $r'_i>r_i$ and a proper Levi subgroup $M_i$ of $G_i$ containing $T(F)$ and such that $U_{x_i,r'_i}U_{M_i,x_i,r_i}\subset H\cap G_i\subset(G_i)_{x_i,r'_i}(M_i)_{x_i,r_i}$.
\end{itemize}

(Note that $r_1,\dots,r_k$ are not related in any way, and can be completely different from each other.)

For every $i$ such that $U_{x_i,r_i}$ is of the first kind, set $M_i=G_i$ and $r'_i=r_i$. Set $M=\prod_{i=1}^kM_i$, and $U_{M,x,(r_1,\dots,r_k)}=\prod_{i=1}^kU_{M_i,x_i,r_i}$; we deduce from above the following generalization of the theorem \ref{cln}:

\begin{theo}\label{cln2}
Assume the conditions of the proposition \ref{nsgf} are satisfied. There exists $x\in\overline{A_P}$ and $r_1,\dots,r_k\in{\mth{R}}^*_+$ such that $H$ satisfies one of the following conditions:
\begin{itemize}
\item $U_{x,(r_1,\dots,r_k)}\subset H\subset G_{x,(r_1,\dots,r_k)}$;
\item there exist $r'_1,\dots,r'_k$ such that $r'_i\geq r_i$ for every $i$ and $r'_i>r'_i$ for at least one $i$, and a proper Levi subgroup $M$ of $G$ containing $T$, such that $U_{x,(r'_1,\dots,r'_k)}U_{M,x,(r_1,\dots,r_k)}\subset H\subset G_{x,(r'_1,\dots,r'_k)}M_{x,(r_1,\dots,r_k)}$.
\end{itemize}
\end{theo}


\begin{thebibliography}{99}
\bibitem{bou} N. Bourbaki, Groupes et alg\`ebres de Lie, chapitre 6: Syst\`emes de racines, Hermann.
\bibitem{bt} F. Bruhat, J. Tits, Groupes r\'eductifs sur un corps local, Publications Math\'ematiques de l'IHES, vol. 41 (1972), pp. 5-251 (I), and vol 60 (1984), pp. 5-184 (II).
\bibitem{chev} C. Chevalley, Sur certains groupes simples, Tohoku Math. Journal, vol 7 (2), 1955, pp. 14-66.
\bibitem{jan} J.C. Jantzen, Representations of Algebraic Groups, Pure and Applied Mathematics, vol. 131, 1987.
\bibitem{mp1} A. Moy, G. Prasad, Unrefined minimal K-types for p-adic groups, Inventiones mathematicae, vol. 116 (1994), pp. 393-408.
\bibitem{mp2} A. Moy, G. Prasad, Jacquet factors and unrefined minimal K-types for p-adic groups, Commentarii Mathematici Helvetici, vol. 71 (1996), pp. 98-121.
\bibitem{pr} G. Prasad, M.S. Raghanathan, Topological central extensions of semi-simple groups over local fields (I and II), Annals of Mathematics, vol. 119 (1984), pp. 143-268.
\bibitem{spr} T.A. Springer, Linear algebraic groups, Progress in Mathematics, Birkh\"auser, vol 9, 2001 (2nd edition).
\bibitem{st} R. Steinberg, Lectures on Chevalley groups, Yale University, 1967.
\end{thebibliography}
\end{document}